\documentclass[notitlepage, 11pt]{article}
\usepackage{graphicx}

\usepackage[usenames, dvipsnames]{color}

\usepackage{bbm}

\usepackage{color}
\usepackage{latexsym}
\usepackage{amsmath, amssymb, amsfonts}
\usepackage[toc]{appendix}
\numberwithin{equation}{section}

\newtheorem{lemma}{Lemma}
\newtheorem{theorem}{Theorem}

\newtheorem{definition}{Definition}

\newtheorem{proposition}{Proposition}

\newtheorem{assumption}{Assumption}

\newtheorem{remark}{Remark}

\newcommand{\beginsec}{
\setcounter{lemma}{0}
\setcounter{theorem}{0}
\setcounter{corollary}{0}
\setcounter{definition}{0}
\setcounter{example}{0}
\setcounter{proposition}{0}
\setcounter{condition}{0}
\setcounter{assumption}{0}
\setcounter{conjecture}{0}
\setcounter{problem}{0}
\setcounter{remark}{0}
}

\newcommand{\noi}{\noindent}

\newcommand{\E}{\mathbb{E}}
\newcommand{\R}{\mathbb{R}}

\newcommand{\N}{\mathbb{N}}

\newcommand{\la}{\lambda}

\newcommand{\eps}{\varepsilon}

\newcommand{\ph}{\varphi}

\newcommand{\al}{\alpha}

\newcommand{\del}{\delta}
\newcommand{\om}{\omega}

\newcommand{\Gam}{\mathnormal{\Gamma}}

\newcommand{\PP}{{\mathbb P}}

\newcommand{\calB}{{\cal B}}
\newcommand{\calC}{{\cal C}}
\newcommand{\calD}{{\cal D}}

\newcommand{\calF}{{\cal F}}

\newcommand{\calK}{{\cal K}}
\newcommand{\calM}{{\cal M}}

\newcommand{\calP}{{\cal P}}

\newcommand{\calU}{{\cal U}}

\newcommand{\calPTL}{{\calP_{T,L}}}

\newcommand{\skp}{\vspace{\baselineskip}}

\newcommand\iy{\infty}

\newcommand{\qed}{\hfill $\Box$}

\newcommand{\lani}{\lambda^{n,i}}
\newcommand{\muni}{\mu^{n,i}}
\newcommand{\Qni}{Q^{n,i}}
\newcommand{\tiQni}{\tilde{Q}^{n,i}}

\newcommand{\aln}{\alpha^{n}}

\newcommand{\alni}{\alpha^{n,i}}

\newcommand{\Nni}{N^{n,i}}
\newcommand{\Mni}{M^{n,i}}

\newcommand{\Ani}{A^{n,i}}
\newcommand{\Dni}{D^{n,i}}

\newcommand{\Jni}{J^{n,i}}

\newcommand{\tQni}{\tilde Q^{n,i}}
\newcommand{\tAni}{\tilde A^{n,i}}
\newcommand{\tDni}{\tilde D^{n,i}}

\newcommand{\tAnj}{\tilde A^{n,j}}
\newcommand{\tDnj}{\tilde D^{n,j}}

\newcommand{\tRni}{\tilde R^{n,i}}
\newcommand{\tbni}{\tilde b^{n,i}}

\newcommand{\talni}{\tilde{\alpha}^{n,i}}

\newcommand{\tYni}{\tilde Y^{n,i}}
\newcommand{\tnun}{\tilde\nu^{n}}

\newcommand{\bnu}{\bar\nu}

\newcommand{\tQn}{\tilde Q^{n}}

\newcommand{\limninf}{\underset{n\to\iy}{\lim\;\inf\;}}
\newcommand{\limnsup}{\underset{n\to\iy}{\lim\;\sup\;}}
\newcommand{\limn}{\lim_{n\to\iy}}

\newcommand{\clc}{{\cal C}}
\newcommand{\clb}{{\cal B}}
\newcommand{\clu}{{\cal U}}
\newcommand{\clr}{{\cal R}}
\newcommand{\clo}{{\cal O}}
\newcommand{\clm}{{\cal M}}
\newcommand{\cld}{{\cal D}}
\newcommand{\cls}{{\cal S}}
\newcommand{\clp}{{\cal P}}
\newcommand{\clf}{{\cal F}}

\newcommand{\cla}{{\cal A}}

\title{Rate Control under Heavy Traffic with Strategic Servers\thanks{This is the final version of the paper. To appear in {\it The Annals of Applied Probability}.}}
\author{Erhan Bayraktar\thanks{Research partially supported by the National Science Foundation (DMS-1613170) and the Susan M.~Smith Professorship.} 
\and
Amarjit Budhiraja\thanks{Research partially supported by the National Science Foundation (DMS-1305120), the Army Research Office (W911NF-14-1-0331) and DARPA (W911NF-15-2-0122).}
\and
Asaf Cohen}

\date{\today}

\begin{document}

\maketitle

\begin{abstract}
\noi We consider a large queueing system that consists of many strategic servers that are weakly interacting. Each server processes jobs from its unique
critically loaded buffer and controls the rate of arrivals and departures associated with its queue to minimize its expected cost. The rates and the cost functions in addition to depending on the control action, can depend, in a symmetric fashion, 
on the size of the individual queue and the empirical measure of the states of all queues in the system. 
In order to determine an approximate Nash equilibrium for this finite player game we construct a Lasry-Lions type mean-field game (MFG) for certain reflected diffusions that governs the limiting behavior. 
Under conditions, we establish the convergence of the Nash-equilibrium value for the finite size queuing system to the value of the MFG. 
\skp

\noi{\bf AMS Classification:} 60K25, 91A13, 60K35, 93E20, 60H30, 
60F17.
\skp

\noi{\bf Keywords:} Heavy traffic limits, queuing systems, strategic servers, mean-field games, rate control, reflected diffusions.
\end{abstract}

\section{Introduction}\label{sec1}
\beginsec

 Rate controlled queueing systems commonly arise from applications in communication systems, see e.g.~\cite{Fendick1994,Kushner2001,Ata2005,Budhiraja2011} and references therein.
Recently, they have also been considered in modeling limit order books, see e.g. \cite{MR2838580,MR3274927,MR3358591,Bowe,Blair,Lachapelle2015}.
A common approach to the study of such rate control problems when the system is in heavy traffic is through diffusion approximations. In a problem setting where there is interaction between servers/queues in that the rates or costs
associated with a particular queue and server can depend on the states of the other queues, this approach leads to a stochastic control problem for $n$-dimensional reflected diffusions, where $n$ is the number of queues in the system. When $n$ is large such control problems are computationally intractable and in general this `curse of dimensionality' is unavoidable. However, when there are certain symmetries present and the interaction between
queues is weak, in that each queue has $\clo(1/n)$ affect on any other queue in the system, a natural approach is to consider, in addition to heavy traffic, another asymptotic regime where the number of queues $n$ approaches $\infty$ as well. Such model settings arise in many applications, e.g., cloud computing, live streaming, limit order books, customer service systems, etc. In many of these contexts the servers are strategic, for example, in customer service networks, servers respond to workload incentives (see \cite{Gopalakrishnan2014}), and in the context of limit order books buyers and sellers place their orders in a strategic manner and interact weakly through their impact on the price distribution.
\skp

\noi We study one of the simplest forms of queuing networks, namely a collection of $n$ single server queues. 
The system is assumed to be \emph{critically loaded}. For this, we fix an arbitrary sequence of scaling parameters $\{e_n\}_n$ that satisfy $e_n\to\iy$ when $n\to\iy$. We assume that, for each queue, the traffic intensity, which is the ratio between the rate of arrivals to the rate of service is $1-\clo(1/\sqrt{e_n})$ (see \cite{Chen2001}). The arrival and service rates for each queue are $\clo(e_n)$. 
In absence of control and state dependence of rates, the analysis simply reduces to that of a single $M/M/1$ queue which in the heavy traffic can be approximated by a one dimensional reflected Brownian motion. 
In the setting we consider, every server can exercise control on the arrival and service rates associated with its own queue. This control which is $\clo(\sqrt{e_n})$ is of lower order compared to the overall rate but it can have significant impact on performance in the asymptotic regime we consider. 
In the heavy traffic regime with a fixed number of servers, performance improvement using an $\clo(\sqrt{e_n})$ control has been well studied (see for example \cite[Chapter 9]{Kushner2001}). In a setting (that is quite different from the one considered here) where the number of nodes/servers approach $\infty$, numerical results that show performance improvement under $\clo(\sqrt{e_n})$ controls can be found in \cite{BudFri2016}.
In the model we consider, the rates can depend on the state of the individual queue, furthermore a particular queue's state is influenced by the remaining queue states through their empirical measure. The control action for each server can use information on the history of queue lengths, arrival and processing times and control actions for all the queues in the system. We consider a rate control problem where each server aims to minimize its individual cost. Although many different types of cost criteria can be considered, here for simplicity we consider a cost function over a finite time horizon. This cost function may depend 
on the individual queue lengths, the control action, and the empirical measure of all the queue states. A natural way to formulate optimality for such $n$-player games is through the notion of a (near) Nash equilibrium. Computing Nash equilibria, even if one considers the simplified heavy traffic approximation in terms of reflected diffusions, is computationally a very challenging problem when $n$ is large. The goal of this work is to analyze an asymptotic formulation where simultaneously each queue approaches heavy traffic {\em and} the number of queues approaches $\infty$. There is extensive literature on heavy traffic limits of rate controlled queuing systems with a fixed number of queues \cite{Kushner2001,Ata2005,Budhiraja2011}, however to the best of our knowledge this is the first work to study the asymptotics for controlled queues where in addition to heavy traffic, the number of queues approach infinity as well. We will show that the asymptotics in this regime is governed by a Lasry-Lions type mean field game (MFG) for reflected diffusions. We will use the solution of the MFG to construct an asymptotic Nash equilibrium for the $n$-player queueing system. This equilibrium has a simple and appealing decentralized structure: each server bases his decision only upon the length of his own queue and a deterministic measure valued function obtained from the solution of the MFG. We also prove that the value 
of the above near Nash equilibrium converges to the value function of the MFG as $n\to \infty$.
In general closed form solutions for MFG of the form that arise here are not available and thus one needs numerical approximations. In \cite{BBC2017} we study one such procedure that uses the Markov chain approximation method (\cite{Kushner1992}) and establish convergence of the scheme over a small time interval. We refer the reader to \cite{Achdou2010,Lachapelle2010, Gueant2012,Achdou2013,  Carlini2014, Carlini2015,Benamou2015, Achdou2016, Chassagneux2017} for some recent results on numerical methods for mean field games.

\skp

\noi The theory of mean field games was initiated a decade ago in the seminal work of Lasry and Lions \cite{Lasry2006,Lasry2006b,Lasry2007}, and Huang, Malham{\'e}, and Caines \cite{Huang2006,Huang2007}. In recent years there has been a growing interest in this field. For recent theoretical developments and applications of this theory see \cite{Cardaliaguet2013,Gueant2013,Gomes2013,Carmona2013,Lacker2015,Carmona2015,Lacker2015general,Fischer2016} and references therein. Mean field approximations for weakly interacting stochastic particles have a long history starting from the works of Boltzmann, McKean, Kac and others (see \cite{sznit} and references therein). Even in the context of queuing systems and communication networks, there have been many works \cite{gibbens1990,Baccelli1992,Hunt1994,Vvedenskaya1997,Borovkov1998,graham2000, Benaim2008,Bobbio2008,graham2009,antunes2008, BudFri2016}.
 Another related branch is agent based models with mean-field interaction (but without strategic agents), see e.g.~\cite{MR1856682,BayraktarHorst1,BayraktarHorst2}. MFGs have also been used by queuing theorists as a tool in recent years see e.g.~\cite{Manjrekar2014,Wiecek2015,Li2015,Lachapelle2015}. In contrast to these works we consider an MFG associated with queuing systems in heavy traffic.
\skp

\noi We now make some comments on proof techniques. Roughly speaking, solution of an MFG considered here can be viewed as the solution of a fixed point problem on the space of probability measures on certain path spaces (see Section \ref{sec3}).
In order to characterize such solutions, in a setting where each agent's state evolution is described through a diffusion
process, Lasry and Lions \cite{Lasry2006,Lasry2006b,Lasry2007} studied wellposedness of two coupled nonlinear partial differential equations; one is an equation of Hamilton-Jacobi-Bellman (HJB) type while the second takes the form of a Kolmogorov forward equation. A somewhat different approach is taken in recent works of Carmona, Delarue, and Lacker \cite{Carmona2013, Carmona2015}. Using probabilistic methods, authors characterize the MFG solution as a solution to certain forward backward stochastic differential equations. 

\noi In all the above papers the $n$-player system is described through a collection of stochastic differential equations and the mean field game gives the asymptotic behavior of the system as $n\to \infty$.
In contrast, in the current work, for a fixed value of $n$ the state process is given through a collection of controlled jump-Markov processes with jump sizes that approach $0$ as $n\to \infty$.
Thus here we have two forms of asymptotic behavior, one corresponds to the large agent limit while the other corresponds to diffusion approximations of pure jump processes with vanishing jump sizes. Analyzing this simultaneous limit behavior, which is key to the proof of Theorem \ref{thm1} that identifies an asymptotic Nash equilibrium for the $n$-player game as $n\to \infty$, requires new techniques. A key ingredient
 is to prove suitable tightness properties and to analyze and characterize the weak limit points of certain stochastic processes and random measures.  Here we make use of a result from \cite{Kotelenez2010} (see Lemma 4.2 therein) which can be viewed as an extension of de Finetti's theorem for sequences of random measures on certain path spaces. In proving suitable tightness properties of control processes we consider a relaxed control formulation. Our assumption on the uniqueness of the minimizer (Assumption \ref{assumption1}(b)) ensures that extending the class of controls in this manner does not lower the cost. Such relaxed control formulations have been widely used in stochastic control theory, see e.g.~\cite{Fleming1976, El1987, Borkar1989, Kushner1992}. More recently, they have also been invoked in the study of mean field games \cite{Lacker2015, Lacker2015general, Fischer2016}.
Other key ingredients needed for the proof are (a) regularity properties of the optimal feedback controls for certain 
 stochastic control problems for reflected diffusions; (b) weak formulation of stochastic control problems for reflected diffusions. We find that the approach based on HJB characterizations of value functions is particularly well suited for our problem setting. In particular, regularity results from \cite{Lieberman1996} for quasilinear PDE with Neumann boundary conditions give us the required estimates for obtaining the desired properties of the optimal feedback controls.
Use of HJB theory in our analysis is a feature that is common to the approach initiated by Lasry and Lions, however a
point of departure is that instead of using Kolmogorov forward equations we characterize probability distributions as unique weak
solutions of suitable reflected diffusions. In this sense our approach is closer to \cite{Carmona2015} although, unlike \cite{Carmona2015}, we do not make use of forward backward stochastic differential equations in our work.
\skp

\noi In summary our main contributions are as follows. We
\begin{itemize} \itemsep0em
\item consider a rate control problem for large symmetric queuing systems in heavy traffic with strategic servers;
\item introduce an MFG for controlled reflected diffusions, and in Theorem \ref{thm_fixed} establish its solvability under Assumption \ref{assumption1}, and prove unique solvability assuming in addition Assumption \ref{assumptionU};
\item use the solution of a diffusion MFG to construct, under Assumptions \ref{assumption1} and \ref{assumption4}, an asymptotically optimal Nash equilibrium for the $n$ player countable state
game using techniques that combine  heavy traffic analysis with large agent asymptotics (Theorem \ref{thm1}).
\end{itemize}
Thus  Assumption \ref{assumption1} is a basic assumption for all our results. This assumption imposes Lipschitz continuity of the various functions in the model. Assumption \ref{assumptionU}, that is introduced for the  uniqueness of the solution of the MFG,  is common in the literature of MFG. It says that the drift is independent of the mean-field term and that the running cost and terminal cost satisfy a certain monotonicity property, see Remark \ref{rem_uniqueness} for a discussion of the assumption. In order to obtain an asymptotic Nash-equilibrium, we require, in addition to Assumption \ref{assumption1},
 Assumption \ref{assumption4} which  as the first part of Assumption \ref{assumptionU} says that the drift is free of the mean-field term and in addition requires  a basic convergence property for the initial conditions. For this result  we do not require the monotonicity  property, in particular the uniqueness of the solution of the MFG is not needed. 

\noi The paper is organized as follows. In Section \ref{sec2} we introduce the queueing model, the scaling regime and the cost criterion. Next in Section \ref{sec3} we introduce the MFG and present our main result on its solvability (Theorem \ref{thm_fixed}). Section \ref{sec4} constructs an asymptotic Nash equilibrium for the $n$-player game and the main result of this section is Theorem \ref{thm1} which proves asymptotic optimality. 

\subsection{Notation} \label{sec11}
We use the following notation.
For metric spaces $\cls_1, \cls_2$ denote by $\clc(\cls_1:\cls_2)$ the space of $\cls_2$ valued continuous functions on $\cls_1$. When $\cls_2 = \R$, we abbreviate the notation to simply $\clc(\cls_1)$. The space $\clc([0,T]: \cls)$ for a Polish space $\cls$ will be equipped with the
 uniform topology. We will denote by $\cld([0,T]: \cls)$ the space of $S$ valued functions that are right continuous and have left limits (RCLL) defined on $[0,T]$. This space is equipped with the usual Skorohod topology.
For $f \in \cld([0,T]:\R^d)$ and $0\le t\le  T$, $\|f\|_t \doteq \sup_{0\le s \le t}\|f(s)\|$. 
In case that $d=1$, we often use $|f|_t$. 
We will denote by $\text{Lip}_1(\cls)$ the space of real Lipschitz functions on $\cls$ whose Lipschitz norm is bounded by $1$, namely the class of functions $f:\cls \to \R$ with
$$|f(x) - f(y)| \le d(x,y),\; x,y \in \cls$$
where $d$ denotes the metric on $\cls$.
Denote by $\calP(\cls)$ the space of probability measures on $\cls$. We endow $\calP(\cls)$ with the topology of weak convergence of measures. Convergence in distribution of $S$ valued random variable $X_n$ to $X$ will be denoted as
$X_n \Rightarrow X$. For $T, L \in (0,\infty)$, the space $\clp(\clc([0,T]:[0,L]))$ will be denoted as $\clp_{T,L}$.
The Wasserstein distance of order $1$ on $\clp(\cls)$, where $\cls$ is a compact metric space,
 is defined as
\begin{align}\notag
W_1(\eta',\eta)=\inf\left\{\left[\int_{\cls}d(x,y)d\pi(x,y)\right] : \pi\in\calP(\cls\times \cls)\;\;\text{ with marginals $\eta'$ and $\eta$}\right\},
\end{align}
 where $\eta,\eta'\in\calP(\cls)$.
 We denote by $\calC^{1,2}([0,T]\times [0,L])$ the space of 
 functions $f:(0,T)\times (0,L)\to \R$ that are continuously differentiable (resp., twice continuously differentiable) with respect to (w.r.t.) the first (resp., second) variable and are such that the function and the derivatives extend continuously to $[0,T]\times [0,L]$. For $\phi \in \calC^{1,2}([0,T]\times [0,L])$, $D_t\phi, D\phi, D^2\phi$ will denote the time derivative and the first two space derivatives of $\phi$, respectively.
For $x \in \cls$, $\delta_x \in \clp(\cls)$ denotes the 
Dirac measure at $x$.

\section{The $n$-server queuing control problem}\label{sec2}\beginsec
We consider a symmetric $n$-server stochastic queueing system. Each server $i \in \{1, \dots, n\}$ 
is associated with a queue with a finite capacity and is able to control the rates of the service and arrivals of jobs to its queue. 
The rates can also depend on the state of the queue and on the empirical measure of the states of all the $n$-queues. We will consider a regime where the arrival and processing rates are approximately the same and are of the same order as the number of queues in the system. We now describe our precise scaling regime and introduce the various processes that determine the evolution of the state of the system.

\subsection{Diffusion scaling}\label{sec2a}
 Fix $T,L>0$ and $n\in\N$. Here $T$ denotes the terminal time of our finite time horizon and $[0,L]$ will be the state space of the scaled queue length process.
Let $U$ be a compact subset of $\R$. We will denote by $\lani$ and $\muni$ the (controlled) arrival and service rates associated with queue $i$. 
Fix an arbitrary sequence of scaling parameters $\{e_n\}_n$ that satisfy $e_n\to\iy$ when $n\to\iy$. The rates will be $\clo(\sqrt{e_n})$ perturbations of certain nominal (uncontrolled) $\clo(e_n)$ overall rates.
More precisely, we assume that there exist $\hat\la,\hat\mu>0$, and bounded and measurable functions $\la,\mu: [0,T]\times\calP([0,L])\times[0,L]\times U \to \R$ such that 
\begin{align}\notag
\lani(t)&=\hat\la e_n+\la(t,\tnun(t), \tiQni(t), \alni(t))\sqrt{e_n}+o(\sqrt{e_n}),\\
\muni(t)&=\hat\mu e_n+\mu(t,\tnun(t),\tiQni(t), \alni(t))\sqrt{e_n}+o(\sqrt{e_n}), \label{2}
\end{align}
where $\tiQni(t) = \tfrac{1}{\sqrt{e_n}} \Qni(t)$, $\Qni(t)$ is the size of the $i$-th queue at time $t$,
$$\tnun(t)\doteq\frac{1}{n}\sum_{i=1}^n\delta_{\tiQni(t)}$$
 is the empirical distribution of the scaled queue lengths in the $n$-th system at time $t$, and $\alni(t)$ is the control that server $i$ exercises at time $t$. 
The term $o(\sqrt{e_n})$ represents an expression of the form $r_n(t,\tnun(t), \tiQni(t), \alni(t))$
where $r_n:[0,T]\times\calP([0,L])\times[0,L]\times U \to \R$ are measurable functions such that
$r_n/\sqrt{e_n}$ converges uniformly to $0$ as $n\to \infty$.
Additional conditions on $\la$ and $\mu$ will be introduced in later sections.
Each server in the $n$-th system has a finite buffer of size
$L^n=\sqrt{e_n}L$ and arriving jobs to a full buffer are lost. We assume that the system is in heavy traffic, namely 
$$\hat\la=\hat\mu.$$

We now give a precise description of controlled stochastic processes of interest.
For each fixed $n$, let $(\Omega',\calF',\PP')$ be a probability space on which are given unit rate independent Poisson processes $\Nni$ and $\Mni$, $i = 1,\ldots,n$. Roughly speaking, $\Nni$ will correspond to the stream of jobs entering the $i$-the queue and $\Mni$ to the jobs that leave the system after service. The evolution of the $i$-th controlled queue is given as follows.
\begin{align}\label{3}
\Qni(t)=\Qni(0)+\Ani(t)-\Dni(t), \quad i = 1,\ldots,n,\quad t\in[0,T],
\end{align}
where $\Qni(0)\ge 0$ is the initial size of the $i$-th queue, 
\begin{equation}
	\label{eq:eq928}\Ani(t)=\Nni\left(\int_0^t 1_{\{\Qni(s)<L\}}\lani(s)ds\right),\; \Dni(t)=\Mni\left(\int_0^t 1_{\{\Qni(s)>0\}}\muni(s)ds\right)\end{equation}
are the arrival and departure processes respectively, where $\lani$ and $\muni$ are defined as in \eqref{2} in terms of control processes $\{\alpha^{n,i}\}$.
%
%
We will assume that $\{\Qni(0)\}_{i=1}^n$ are exchangeable for all $n$. We require that the control processes are suitably non-anticipative. Specifically, we assume that for some filtration $\{\calF_t\}$ on $(\Omega',\calF',\PP')$, $\alni$ is $\{\calF_t\}$-progressively measurable, $\Qni$, $n\ge i\ge 1$ are $\{\calF_t\}$-adapted, and 
\begin{align}\notag
\tAni(t)\doteq\frac{\Ani(t)-\int_0^t 1_{\{\Qni(s)<L\}}\lani(s)ds}{\sqrt{e_n}},\; 
\tDni(t)\doteq\frac{\Dni(t)-\int_0^t 1_{\{\Qni(s)>0\}}\muni(s)ds}{\sqrt{e_n}}
\end{align}
 are $\{\calF_t\}$ martingales with quadratic variations
\begin{equation}\label{eq: eq415}
\langle\tAni, \tAnj \rangle(t) =\delta_{ij}\frac{1}{e_n}\int_0^t 1_{\{\Qni(s)<L\}}\lani(s)ds,\; \langle \tDni, \tDnj \rangle(t) = 
\delta_{ij}\frac{1}{e_n}\int_0^t 1_{\{\Qni(s)>0\}}\muni(s)ds,
\end{equation}
and $\langle \tAni, \tDnj \rangle(t) = 0$, $t\in [0,T]$, $i,j = 1, \ldots, n$ where $\delta_{ij}=1$ if $i=j$ and $0$ otherwise.
The process $\aln = \{\alni\}_{i=1}^n$ will be referred to as an {\it admissible} control and we denote the collection of all such controls by ${\calU}^n$. With an abuse of terminology, for $\aln = \{\alni\}_{i=1}^n$ as above, we will refer to
$\alni$ as an admissible control for the $i$-th player.

%
From \eqref{3} we have the following evolution equation for the scaled queue length processes. For $t\in [0,T]$
\begin{align}\label{7}
\tQni(t)=\tQni(0)+\tAni(t)-\tDni(t)+\int_0^t\tbni(s)ds+\tYni(t)-\tRni(t) + o(1),
\end{align}
where 
\begin{align}\label{8}
\tYni(t)\doteq\frac{\int_0^t 1_{\{\tQni(s)=0\}}\muni(s)ds}{\sqrt{e_n}},\;
\tRni(t)\doteq\frac{\int_0^t 1_{\{\tQni(s)=L\}}\lani(s)ds}{\sqrt{e_n}},
\end{align}
\begin{align}\notag
\tbni(t)\doteq b(t,\tnun(t),\tQni(t), \alni(t)), \quad b\doteq\lambda-\mu,
\end{align}
and $o(1)$ represents a stochastic process $\eta^{n,i}$ satisfying for each $i$
$\sup_{0\le t \le T}|\eta^{n,i}(t)|\to 0$, in probability, as $n \to \infty$.
The above dynamics can equivalently be described in terms of a Skorohod map as we do below.
Let $\cld_0[0,T]$ be the subset of $\cld([0,T]:\R)$ consisting of all $\psi$ such that $\psi(0) \in [0,L]$.
\begin{definition}\label{def_Skorohod}
Given $\psi\in\calD_0[0,T]$ we say the triplet of functions $(\ph,\zeta_1,\zeta_2)\in\calD([0,T]: \R^3)$ solves the Skorohod problem for $\psi$ if the following properties are satisfied:

\noi
(i) For every $t\in[0,T], \;\ph(t)=\psi(t)+\zeta_1(t)-\zeta_2(t)\in [0,L]$.

\noi
(ii) $\zeta_i$ are nonnegative and nondecreasing, $\zeta_1(0)=\zeta_2(0)=0$, and 
\begin{align}\notag
\int_{[0,T]}1_{(0,L]}(\ph(s))d\zeta_1(s)=\int_{[0,T]}1_{[0,L)}(\ph(s))d\zeta_2(s)=0.
\end{align}
We denote by $\Gamma(\psi)=(\Gamma_1,\Gamma_2,\Gamma_3)(\psi)=(\ph,\zeta_1,\zeta_2)$ and refer to $\Gamma$ as the Skorohod map. 
\end{definition}
It is known that there is a unique solution to the Skorohod problem for every $\psi \in\calD_0[0,T]$ and so the Skorohod map in Definition \ref{def_Skorohod} is well defined.
The Skorohod map has the following Lipschitz property (see \cite{Kruk2007}).
\begin{lemma}\label{lem_Skorohod}
There exists $c_S \in (0,\infty)$ such that
for all $\om,\tilde\om\in\calD_0([0, T])$,
\begin{equation}\notag
\sum_{i=1}^3\|\Gam_i(\om) - \Gam_i(\tilde\om)\|_{T}
\le c_S\|\om-\tilde\om\|_T.
\end{equation}
\end{lemma}
The dynamics in \eqref{7} can be described in terms of the Skorohod map as follows
\begin{align}\label{eq_23}
(\tQni,\tYni,\tRni)(t)=\Gamma\left(\tQni(0)+\int_0^\cdot\tbni(s)ds+\tAni(\cdot)-\tDni(\cdot) + o(1)\right)(t),\quad t\in[0,T].
\end{align}
\subsection{The control problem}
\label{sec:sec-cont}
We now introduce the cost function and the control problem studied in this work.
The total expected cost for server $i$ associated with the initial condition $\\$ $\tilde{Q}^n(0)=(\tilde{Q}^{n,1}(0),\ldots,\tilde{Q}^{n,n}(0))$ and control $\aln=(\al^{n,1},\ldots,\al^{n,n})\in \calU^n$ is given by
\begin{align}\label{12}
\Jni(\tQn(0);\aln) &\doteq 
\E\Big[\int_0^T f(t,\tnun(t), \tQni(t), \alni(t))dt + g(\tnun(T),\tQni(T))\\\notag
&\quad+\int_0^Ty(t,\tnun(t))d\tYni(t)+\int_0^Tr(t,\tnun(t))d\tRni(t)\Big],
\end{align}
where $f:[0,T]\times \calP([0,L])\times[0,L] \times U\to\R$ is the running cost, $g: \calP(\R)\times [0,L]\to\R$ is the terminal cost, and $r,y:[0,T]\times \calP([0,L])\to\R_+$ are the costs associated with rejection
of jobs and empty buffers, respectively. Here $f,g,y$, and $r$ are bounded measurable functions that will be required to satisfy additional conditions that will be introduced in Section \ref{sec3} (see Assumption \ref{assumption1}). 
Each player (server) seeks to minimize its cost. A natural formulation of optimality for such $n$-player games is given in terms of a (near) Nash equilibrium. Computing near Nash equilibria for such complex and large multi-player games is
in general intractable and thus we instead consider an asymptotic formulation of the problem.
\begin{definition}
A sequence of admissible controls $\{\talni: 1\le i\le n\}_{n \in \N}$ is called an {\it asymptotic Nash equilibrium} if for every player $j$, and every sequence of admissible controls $\{\beta^n\}_{n=1}^\iy$ for that player, one has
\begin{align}\notag
\limnsup J^{n,j}(\tQn(0);\tilde \al^{n,1},\dots,\tilde \al^{n,n})\le\limninf J^{n,j}(\tQn(0);\tilde \al^{n,1},\dots,\tilde\al^{n,j-1},\beta^n,\tilde\al^{n,j+1},\ldots,\tilde \al^{n,n}).
\end{align}
\end{definition}
Objective of this work is to show that, under conditions, an asymptotic Nash equilibrium exists which can be approximated by analyzing a related MFG.
The main results are Theorem \ref{thm_fixed} (solvability of MFG) and Theorem \ref{thm1} (asymptotic optimality).

\section{The MFG}\label{sec3} \beginsec
A natural approach for constructing asymptotic near Nash equilibria for the above $n$-player game has emerged from
the works of \cite{Cardaliaguet2013,Carmona2013,Carmona2015}. Starting point in this approach is to consider an MFG that formally corresponds to the limit of the above $n$-player games as $n\to \infty$. In this section we give a precise description of this MFG in the current context and give our main results on existence and uniqueness of solutions.
\subsection{Description of the MFG}\label{sec3a} 
 The basic idea in the formulation of the MFG is to approximate the scaled queue length process for a typical queue by a suitable drift-controlled reflected Brownian motion. We next introduce this controlled process.
Let $(\Omega,\calF,\{\calF_t\},\PP)$ be a filtered probability space on which is given a one dimensional standard $\calF_t$-Brownian motion $B$. We will refer to the collection $(\Omega,\calF,\{\calF_t\},\PP, B)$
as a system and denote it by $\Xi$. 
Given $x \in [0,L]$, $t \in [0,T]$, and $\nu \in \clp_{T,L}$, we denote by $\cla(\Xi, t,x,\nu)$ the collection of
all pairs $(\alpha, Z)$ where $\alpha = \{\alpha(s)\}_{0\le s \le T-t}$ is a $U$-valued
$\clf_s$-progressively measurable process, $Z = \{Z(s)\}_{0\le s \le T-t}$ is a $[0,L]\times \R_+\times \R_+$ valued $\clf_s$-adapted continuous process such that, $Z=(X,Y,R)$ and 
\begin{align}\label{13}
Z(s) = (X,Y,R)(s)=\Gamma\left(x+\int_0^\cdot\bar b(u)du+\sigma B(\cdot)\right)(s),\quad s\in[0,T-t],
\end{align}
where 
\begin{align}\notag
\bar b(u)\doteq b(t+u,\nu(t+u),X(u), \al(u)),\; u \in [0, T-t],
\end{align}
$\nu(s)$ is the marginal of $\nu$ at time instant $s$ and $\sigma=\sqrt{2\hat\lambda}$.
We now introduce the cost function in the MFG. Given $\nu \in \clp_{T,L}$, $x \in [0,L]$, $t \in [0,T]$, and a system
$\Xi$ as above, let $(\alpha,Z) \in \cla(\Xi, t,x,\nu)$. Define
 \begin{align}\notag
J_\nu(t,x,\al,Z)&\doteq \E\Big[\int_0^{T-t} f(s+t,\nu(s+t),X(s),\al(s) )ds + g(\nu(T),X(T-t))\\
&\quad+\int_0^{T-t}y(s+t,\nu(s+t))dY(s)+\int_0^{T-t}r(s+t,\nu(s+t))dR(s)\Big].\label{eq:eq234}
\end{align}
We define the value function associated with the above cost as:
\begin{align}\label{15}
V_\nu(t,x)=\inf_{\Xi}\inf_{(\al,Z)\in \cla(\Xi,t,x,\nu)}J_\nu(t,x,\al,Z).
\end{align}
\begin{definition}
	\label{def:mfg729}
A {\bf solution to the MFG} with initial condition $x \in [0,L]$ is defined to be a $\nu \in \clp_{T,L}$ such that
there exist a system $\Xi$ and an $(\alpha,Z) \in \cla(\Xi,0,x,\nu)$ such that $Z=(X,Y,R)$ 
satisfies $\PP \circ X^{-1} = \nu$ and
\begin{equation}\label{eq:eq716}
	V_\nu(0,x) = J_{\nu}(0, x, \alpha,Z).
\end{equation}
If there exists a unique such $\nu$, we refer to $V_\nu(0,x)$ as the {\bf value of the MFG} with initial condition
$x$.
\end{definition}
To find a solution of the MFG one usually follows the following iterative procedure: 

\begin{itemize}
\item[(i)] For a fixed $\nu \in \clp_{T,L}$ solve the stochastic control problem \eqref{15}(with $t=0$), namely find a system $\Xi$ and $(\al,Z)\in \cla(\Xi,0,x,\nu)$
such that \eqref{eq:eq716} holds. Denote by $\bar \nu$ the law of $X$ where $Z=(X,Y,R)$ and write
$\bar \nu = \Phi(\nu)$ (this is not precise since in general there may be more than one solution of the stochastic control problem in \eqref{15}).
\item[(ii)] Find the fixed point of the map $\Phi$, namely a $\bar \nu \in \clp_{T,L}$ for which
$\bar \nu = \Phi(\bar \nu)$. Note that by definition, such a $\bar \nu$ will be a solution of the MFG.

\end{itemize}

We now analyze the MFG by following the above steps. The main result is Theorem \ref{thm_fixed} which gives existence of solutions of the MFG under suitable conditions and proves uniqueness
of solutions under stronger conditions.

\subsection{Solving the stochastic control problem \eqref{15}}\label{sec3b1} 
For $c \in(0,\infty)$, let $\calM_c$ be the collection of all $\nu \in \clp_{T,L}$
 such that 
\begin{align}\notag
\sup_{0\le s<t\le T}\frac{W_1(\nu(t),\nu(s))}{(t-s)^{1/2}}\le c
\end{align}
and let 
\begin{align}\notag
\clm_0 = \cup_{c>0} \clm_c.
\end{align}
Fix a measure $\nu\in\calM_0$.
For $\nu \in \clp_{T,L}$, the function $t \mapsto \nu(t)$ is a continuous function from $[0,T]$ to
$\clp([0,L])$ and with an abuse of notation we denote this continuous function once more as $\nu$.
As one might expect, the value function $V_{\nu}(t,x)$ corresponds to the Hamilton-Jacobi-Bellman (HJB) equation:
\begin{align}\label{HJB1} 
-D_t \phi-H(t,\nu(t),x,D \phi)-\frac{1}{2}\sigma^2 D^2\phi=0,\qquad (t,x)\in[0,T]\times[0,L],
\end{align}
with the boundary conditions (BC)
\begin{align}\label{HJB2}
\phi(T,x)=g(\nu(T),x),\; D\phi(t,0)=-y(t,\nu(t)),\text{ and } D\phi(t,L)=r(t,\nu(t)),\; t\in[0,T],
\end{align}
where $H$ is the Hamiltonian given as
$$H(t,\eta,x,p)=\inf_{u\in U} h (t,\eta,x,u,p)$$ 
and $h:[0,T]\times\calP([0,L])\times [0,L]\times U\times \R \to \R$ is defined as
\begin{align} \label{19}
h(t,\eta,x,u,p)=f(t,\eta,x,u)+b(t,\eta,x,u)p.
\end{align}
We now introduce a key condition under which the above HJB equation characterizes the value function $V_{\nu}(t,x)$. 
\begin{assumption}\label{assumption1}$\,$
\begin{enumerate}
\item[(a)] There exists $c_L \in (0,\infty)$ such that for every $t,t'\in[0,T]$, $\eta,\eta'\in\calP([0,L])$, $x,x'\in[0,L]$, and $\al,\al'\in U$,
\begin{align}
&|f(t,\eta,x,\al)-f(t',\eta',x',\al')| + |g(\eta,x) - g(\eta',x')|+|b(t,\eta,x,\al)-b(t',\eta',x',\al')|\nonumber\\
&\quad +|y(t,\eta)-y(t',\eta')|+|r(t,\eta)-r(t',\eta')|\nonumber\\
&\quad\quad\le c_L (|t-t'|+W_1(\eta,\eta')+|x-x'|+|\al-\al'|).\label{19b}
\end{align}
\item[(b)]
For every $(t,\eta, x,p) \in [0,T]\times \clp([0,L])\times [0,L]\times \R$, there is a unique
$\hat \al(t,\eta,x,p) \in U$ such that
\begin{align}\label{19a}
\hat \al(t,\eta,x,p)=\underset{u\in U}{\arg\min}\;h(t,\eta,x,u,p) .
\end{align}
\end{enumerate}
\end{assumption}
From Berge's maximum theorem (see \cite[Theorem 17.31]{Aliprantis2006}) and part (b) of the above assumption
it follows that
$\hat \al$ is a continuous function on $[0,T]\times \clp([0,L])\times [0,L]\times \R$.
Also note that \eqref{19b} implies that $b,f,g,y,r$ are bounded functions, in particular,
\begin{equation}\sup_{(\eta,x,u)\in [0,T]\times \clp([0,L])\times [0,L]\times U}|b(t,\eta,x,u)| \doteq c_B < \infty.
	\label{eq:eq429}
\end{equation}
The first part of the condition implies H\"{o}lder-1/2 continuity in $t$ when $\eta$ is replaced by $\nu(t)$
for some $\nu \in \clm_0$.
For example, if $\nu \in \clm_c$, for $t,t' \in [0,T]$
\begin{align}\label{eq_20}
|y(t,\nu(t))-y(t',\nu(t'))|&\le c_L(c+\sqrt{T})|t-t'|^{1/2}.
\end{align}
Similar estimates hold for $b,f,r$ and $g$.

Assumption \ref{assumption1} will guarantee that the value function is the unique classical solution of the HJB equation \eqref{HJB1} with BC \eqref{HJB2} (see Proposition \ref{lem2}). 
The assumption will also be used in the fixed point argument (Section \ref{sec3b2}) and in the asymptotic analysis of the $n$-player game (Section \ref{sec4}). 
Using Girsanov's theorem it is easily checked that given a measurable function 
$\gamma: [0,T]\times [0,L] \to U$, $\nu \in \clp_{T,L}$, and $(t,x) \in [0,T]\times [0,L]$
there is a unique weak solution $Z$ to 
\begin{equation}
	\label{eq:feedback}
	Z(s) = (X,Y,R)(s)=\Gamma\left(x+\int_0^\cdot b(t+t', \nu(t+t'), X(t'), \gamma(t', X(t')))dt'+\sigma B(\cdot)\right)(s), 
\end{equation}
$s\in[0,T-t]$, namely there is a system $\Xi = (\Omega,\calF,\{\calF_t\},\PP, B)$ on which is given a $\clf_s$-adapted
continuous process $Z = (X,Y,R)$ solving the above equation and if $\Xi' = (\Omega',\calF',\{\calF'_t\},\PP', B')$ is another system on which is given a $\clf'_s$-adapted
continuous process $Z' = (X',Y',R')$ solving the above equation with $(Z,X, B)$ replaced by
$(Z',X', B')$ then $\PP\circ Z^{-1} = \PP'\circ (Z')^{-1}$. Note also that with $\alpha(s) = \gamma(s, X(s))$,
$(\al,Z)\in \cla(\Xi,t,x,\nu)$. We refer to the function $\gamma$ as a feedback control and we call it an {\em optimal} feedback control for the stochastic control problem in \eqref{15} if
$$V_{\nu}(t,x) = J_{\nu}(t,x,\alpha, Z).$$
The following result which is a consequence of Theorem 13.24 of \cite{Lieberman1996} says that
the HJB equation \eqref{HJB1} with BC \eqref{HJB2} admits a unique classical solution which can be characterized as the value function $V_\nu$. Moreover, there exists an optimal feedback control with certain regularity properties.
We follow the notation from \cite{Lieberman1996}. Let for $\delta \in (0,1]$, $H_{\delta}$ be the collection of
maps $\psi:(0,T)\times (0,L) \to \R$ such that
$$\sup_{0<t <t'<T, 0< x<x'<L} \frac{|\psi(t,x) - \psi(t',x')|}{|t-t'|^{\delta/2}+ |x-x'|^{\delta}} < \infty.$$
Note that such a function can be continuously extended to $[0,T]\times [0,L]$ and we denote the extension by the same symbol.
Also, let $H_{2+\frac{1}{2}}$ be the collection of continuous real functions $\psi$ on $[0,T]\times [0,L]$ 
such that
$x\mapsto \psi(t,x)$ is twice continuously differentiable on $(0,L)$ for all $t \in (0,T)$, 
$t\mapsto \psi(t,x)$ is continuously differentiable on $(0,T)$ for all $x \in (0,L)$,
the functions 
$$\psi,\; D_t \psi,\; D \psi\; D^2 \psi,$$ 
are bounded on $(0,T)\times (0,L)$ and the functions
$D_t \psi$, $D^2 \psi$ are in $H_{1/2}$.
\begin{proposition}\label{lem2}
Fix  $\nu\in \calM_0$ and suppose that Assumption \ref{assumption1} holds. Then $V_\nu\in H_{2+\frac{1}{2}}$ and
it is the unique solution of \eqref{HJB1}--\eqref{HJB2}. 
Furthermore, with $\hat \alpha$ as introduced in Assumption \ref{assumption1}, 
%
the map $(s,x')\mapsto \hat \al(s,\nu(s),x',D V_\nu(s,x'))$ is continuous and the
feedback control $\hat\gamma(u, x')\doteq \hat \al(u+t,\nu(u+t),x', D V_\nu(u+t,x'))$ is an optimal feedback control for \eqref{15} for every $t \in (0,T)$. Moreover, any optimal control $\al$ for \eqref{15} satisfies $\al(u,\omega)=\hat\gamma(u,X(u,\omega))$, $\la_T^t\otimes\PP$ almost surely (a.s.), where $\lambda_T^t$ denotes the Lebesgue measure on $[0,T-t]$. 
\end{proposition}
{\bf Proof:} 
From \cite[Theorem 13.24]{Lieberman1996} and the paragraph following its statement it follows that \eqref{HJB1}--\eqref{HJB2} admits a solution in $H_{2+\tfrac{1}{2}}$. We remark that the key conditions needed to appeal to this theorem are that
$$\sup_{0< t <t' <T} \frac{|y(t,\nu(t)) - y(t',\nu(t'))| + |r(t,\nu(t)) - r(t',\nu(t'))|}{|t-t'|^{1/2}} < \infty,$$
and for each $M\in (0,\iy)$
$$
\sup_{\substack{t,t' \in [0,T],\, x,x' \in [0,L],\, p,p' \in [-M,M],\\ t\neq t', x\neq x', p\neq p'}}
\frac{|H(t,\nu(t), x,p) - H(t',\nu(t'), x', p')| }{|t-t'|^{1/2} + |x-x'| + |p-p'|} < \infty.$$
Both these conditions are easily seen to hold on using Assumption \ref{assumption1} and the property that
$\nu \in \clm_c$ for some $c<\infty$ (see e.g.~\eqref{eq_20}).

We now argue uniqueness and characterize the solution as the value function $V_{\nu}$ in \eqref{15}.
 Let $\psi\in H_{2+ \frac{1}{2}}$ be a solution of \eqref{HJB1}--\eqref{HJB2}. 
Fix $(t,x) \in [0,T]\times [0,L]$, a system $\Xi$ and $(\alpha, Z) \in \cla(\Xi, t,x,\nu)$.
By an application of It\^{o}'s lemma and using \eqref{HJB1}--\eqref{HJB2} we get
\begin{align}
\E[g(\nu(T),X(T-t))] &= \E[\psi(T,X(T-t))]\nonumber\\
& =\psi(t,x) + \E\left[\int_0^{T-t} D\psi(s+t,0)dY(s)-\int_0^{T-t} D\psi(s+t,L)dR(s)\right]\nonumber\\
&\quad +\E\int_0^{T-t}\Big[D_s \psi(s+t,X(s))+b(s+t,\nu(s+t),X(s),\al(s))D \psi(s+t,X(s))\nonumber\\
&\quad +\frac{1}{2}\sigma^2D^2\psi(s+t,X(s))\Big]ds \nonumber\\
&\ge \psi(t,x)+\E\Big[-\int_0^{T-t} f(s+t,\nu(s+t),X(s),\al(s))ds\nonumber\\
&\quad-\int_0^{T-t}y(s+t,\nu(s+t))dY(s) -\int_0^{T-t}r(s+t,\nu(s+t))dR(s)\Big].
\label{eq:eq624}
\end{align}
Hence, with $Z=(X,Y,R)$, 
\begin{align}\label{unique1}
J_\nu(t,x,\al, Z)&=\E\left[\int_0^{T-t}f(s+t,\nu(s+t),X(s),\al(s))ds +g(\nu(T),X(T-t)) \right.\\\notag
&\qquad\qquad\left.
+\int_0^{T-t}y(s+t,\nu(s+t))dY(s) +\int_0^{T-t}r(s+t,\nu(s+t))dR(s)\right]\\\notag
&\ge \psi(t,x).
\end{align}
Since $\Xi, \al, Z$ are arbitrary, we get that 
$V_\nu(t,x)\ge\psi(t,x)$ for all $(t,x)\in [0,T]\times [0,L]$. 
Let $\hat \gamma$ be as in the statement of the proposition with $V_{\nu}$ replaced by $\psi$. 
Let $Z = (X,Y,R)$ be a solution of \eqref{eq:feedback} with $\gamma$ replaced by $\hat \gamma$ given on some system
$\Xi$.
Then a similar calculation using It\^{o}'s formula but with $\alpha(u) = \hat\gamma(u, X(u))$ shows that
$J_{\nu}(t,x,\alpha, Z) = \psi(t,x)$ for all $(t,x) \in [0,T]\times [0,L]$. This shows that $\psi = V_{\nu}$
and so $V_\nu$ is the unique $H_{2+\frac{1}{2}}$ solution of \eqref{HJB1}--\eqref{HJB2}. Also, we have that $\hat \gamma$ 
as in the statement of the lemma is an optimal feedback control. 

Since $V_{\nu}\in H_{2 + \frac{1}{2}}$, the continuity of the map $(t,x)\mapsto \hat \al(t,\nu(t),x,D V_\nu(t,x))$ 
is immediate from the continuity of the map $\hat\al$ that was noted below Assumption \ref{assumption1}.

We now show that $\hat\gamma$ is the unique optimal control for \eqref{15}. Fix an optimal control $(\al, Z) \in \cla(\Xi, t, x, \nu)$ given on some system $\Xi$. We claim that 
$\al(u,\omega)=\hat\gamma(u, X(u,\omega))$, $\lambda_T^t\times \PP$-a.s. Indeed, suppose that there is a set with a positive $\lambda_T^t\times\PP$-measure on which 
the equality fails. Then by Assumption \ref{assumption1}.(b) together with \eqref{HJB1}--\eqref{HJB2}, it follows that \eqref{eq:eq624} holds with a strict inequality, which in turn implies that \eqref{unique1} holds with a strict inequality. Recalling that $V_\nu=\psi$, we arrive at a contradiction.
\hfill$\Box$

\subsection{Solving the MFG}\label{sec3b2}
We now turn to step (ii) in solving the fixed point problem. Although for existence of a fixed point Assumption \ref{assumption1} will suffice, in order to get uniqueness, we will need the following additional assumption. 
 Similar assumption has been used to argue uniqueness of fixed points in previous works on MFG (see e.g.~\cite[Theorem 2.4]{Lasry2007}, \cite[Section 3]{Cardaliaguet2013}, \cite[equation (17)]{Gomes2013}, \cite[Assumption (U)]{Carmona2015}).
Fix $\eta_0 \in \clp([0,L])$.
\begin{assumption}\label{assumptionU}
For every $(t,\eta,x,u)\in[0,T]\times\calP([0,L])\times[0,L]\times U$, 
\begin{align}\label{19q}
b(t,\eta,x,u)&=b(t,\eta_0,x,u),\quad
f(t,\eta,x,u)=f_0(t,\eta,x)+f_1(t,x,u),\\\label{19pp}
y(t,\eta)&=y(t,\eta_0),\quad
r(t,\eta)=r(t,\eta_0).
\end{align}
Moreover, for every $t\in[0,T]$ and $\eta,\eta'\in\calP([0,L])$, $f_0$ and $g$ satisfy the following monotonicity property
 \begin{align}\notag
 \int_0^L [f_0(t,\eta,x)-f_0(t,\eta',x))]d(\eta-\eta')(x)&\ge 0,\\\notag
 \int_0^L (g(\eta,x)-g(\eta',x))d(\eta-\eta')(x)&\ge 0.
 \end{align}
\end{assumption}
Abusing notation, when Assumption \ref{assumptionU} holds, we will write $b(t,x,u) =b(t,\eta_0,x,u)$, $y(t)=y(t,\eta_0)$, and $r(t)=r(t,\eta_0)$.

\begin{remark}\label{rem_uniqueness} Examples satisfying Assumption \ref{assumptionU} are given in \cite[page 8]{Cardaliaguet2013} and \cite[page  6]{bayzha2016}. Another natural example  that satisfies this assumption is a cost function that is linear in the mean-field term. That is,
\begin{align}\notag
f_0(t,\eta,x)=a_1(t)(c_1+\psi_1(x))\int_0^L\psi_1(y)d\eta(y),\qquad g(\eta,x)=(c_2+\psi_2(x))\int_0^L\psi_2(y)d\eta(y),
\end{align}
where $a_1: [0,T]\to\R_+$, $\psi_1,\psi_2:[0,L]\to\R$, and $c_1,c_2\in\R$. 
From a modeling perspective, by choosing positive and nondecreasing $\psi_i$'s and a positive $a_1$, the system planner penalizes all servers collectively when the empirical measure has high 
$\psi_i$-moments and in addition it penalizes  individual servers for long queues. 
\end{remark}

\begin{theorem}\label{thm_fixed}
Under Assumption \ref{assumption1}, there exists a solution of the MFG. If in addition Assumption \ref{assumptionU} holds then there is a unique MFG solution.
\end{theorem}
 The proof will appeal to Schauder's fixed point theorem. Since Schauder's original work (cf.~\cite{Schauder1930}), there have been several versions of this result. We now quote the version (\cite[Theorem 4.1.1]{Smart1974}) that will be used here.
\begin{lemma}\label{lem_Schauder}
Let $\calM$ be a non-empty convex subset of a normed space $\calB$. Let $\Phi$ be a continuous mapping of $\calM$ into a compact set $\calK\subset\calM$. Then $\Phi$ has a fixed point, namely there exists $x \in \calK$ such that
$\Phi(x)=x$.
\end{lemma}

\noi {\bf Proof of Theorem \ref{thm_fixed}:}
We will apply Lemma \ref{lem_Schauder} to the space $\calB$ of finite signed measures on $\calC([0,T]:[0,L])$ which is equipped with the Kantorovich-Rubinstein norm 
\begin{align}\notag
\|\nu\|_{KR}=\sup\left\{\left|\int_{\calC([0,T]:[0,L])}F(\omega)d\nu(\omega)\right|\;:\;F\in\text{Lip}_1(\calC([0,T]:[0,L]))\right\},\; \nu\in\calB .
\end{align}
 The distance driven by the norm, restricted to $\calPTL$, coincides with the Wasserstein's distance of order 1 (see \cite[Remark 6.5]{Villani2009}), which due to compactness of $[0,L]$ generates the same topology
on $\calPTL$ as that for weak convergence.

We now introduce a mapping $\Phi$ on the non-empty and convex set $\calM_0$ that satisfies the conditions stated in Lemma \ref{lem_Schauder}.
%

\noi {\bf Definition of $\Phi$.} 
For $\nu\in \calM_0$, let $\alpha_{\nu}$ denote the optimal feedback control $\hat \gamma$ for \eqref{15} (with $t=0$) given through Proposition \ref{lem2}.
Let $Z^{\nu} = (X^{\nu}, Y^{\nu}, R^{\nu})$ denote the unique weak solution of \eqref{eq:feedback} with $\gamma$ replaced with $\alpha_{\nu}$ given on some system
$\Xi =(\Omega,\calF,\{\calF_t\},\PP, B)$. Define
$\Phi(\nu) = \PP \circ (X^{\nu})^{-1}$. 

\noi {\bf Invariance of $\calM_0$.}
For $\nu\in \calM_0$ and $0\le s \le t <T$
\begin{align}\label{50}
W_1(\Phi(\nu)(t),\Phi(\nu)(s))&\le \E|X^{\nu}(t)-X^{\nu}(s)|\\\notag&\le c_S\max\{c_B,\sigma\}\left((t-s)+\E \sup_{s \le u \le t}|B(u)-B(s)|\right)\\\notag&\le 4\hat{C} (t-s)^{1/2},
\end{align}
where $\hat C \doteq c_S\max\{c_B,\sigma\}$, and the last inequality uses Doob's maximal inequality. This shows that $\Phi(\nu) \in \calM_0$ for all $\nu\in \calM_0$. 

We now show
that $\Phi(\calM_0)$ is contained in a compact set in $\calPTL$, i.e.~$\Phi(\calM_0)$ is relatively compact.

 %
 %
 
\noi {\bf Relative compactness of $\Phi(\calM_0)$.}  For $f \in \clc([0,T]:\R)$ and $\delta >0$, let 
$$\varpi_f(\delta) = \sup_{s,t \in [0,T]:|t-s|\le \delta}|f(t)-f(s)|.$$
Then for $\eps, \delta,\varrho >0$ and $\nu \in \clm_0$, similar to the estimate in \eqref{50},
\begin{align*}
	\Phi(\nu) (X^{\nu}: \varpi_{X^{\nu}}(\delta) \ge \varrho) &\le \frac{1}{\varrho}\E \sup_{0\le s \le t \le s+\delta \le T} |X^{\nu}(t) - X^{\nu}(s)|\\
	&\le \frac{\hat C}{\varrho} \left(\delta + \E \sup_{0\le s \le t \le s+\delta \le T} \sup_{s\le u \le t}|B(u)-B(s)| \right)\\
	&\le \frac{\hat C}{\varrho} (\delta + \E \varpi_{B}(\delta)).
\end{align*}
Since $\E \varpi_{B}(\delta) \to 0$ as $\delta \to 0$, we have from \cite[Theorem 7.3]{Billingsley1999} that $\{\Phi(\nu): \nu \in \clm_0\}$ is relatively
compact in $\calPTL$.

\noi {\bf Continuity of $\Phi$.} 
We now argue that $\Phi$ is a continuous map on $\clm_0$. Consider a system $\Xi =(\Omega,\calF,\{\calF_t\},\PP, B)$ and let $Z^0 = (X^0, Y^0, R^0)$ be given as
\begin{equation}\label{eq:eq433}
	Z^0(s) = (X^0(s), Y^0(s), R^0(s)) = \Gamma(x + \sigma B(\cdot))(s), \; 0 \le s \le T.\end{equation}
Define for $\nu' \in \clm_0$
\begin{align}\label{up1}
\alpha^0_{\nu'}(t) = \hat \al(t,\nu'(t),X^0(t),D V_{\nu'}(t,X^0(t))), \; t \in [0,T].
\end{align}
%
%
%
%
%
%
%
%
Let $\nu^n, \nu \in \calM_0$ be such that $\nu^n \to \nu$. 
Since  $\sup_{0\le t\le T}W_1(\nu^n(t),\nu(t))\le W_1(\nu^n,\nu)$, the above convergence implies that
\begin{align}\label{eq_21}
\limn \sup_{0\le t\le T}W_1(\nu^n(t),\nu(t))=0.
\end{align}
We now show that $\al_{\nu^n}^0\to\al_\nu^0$ in $\lambda_T\times\cal\PP$-measure. 
Recall that $\hat \alpha$ is a continuous map. 
Hence, in view of \eqref{eq_21}, for the desired convergence, it is sufficient to show that 
\begin{align}\label{new79}
\limn\E\left[\int_0^T|D V_{\nu^n}(t,X^0(t)) - D V_{\nu}(t,X^0(t))|dt\right]=0.
\end{align}
The proof is a modification of the proof of Theorem 2.1 in \cite{Hu1997}. 
Applying It\^{o}'s lemma to $V_{\nu'}(t,X^0(t))$ for fixed $\nu'\in\calM_0$ and using from Proposition \ref{lem2} the fact that $V_{\nu'}$ solves \eqref{HJB1}--\eqref{HJB2}, we have for every $t\in [0,T]$, 
\begin{align}\label{eq_40}
&V_{\nu'}(t,X^0(t)) - V_{\nu'}(T,X^0(T)) 
\\\notag
&\quad= \int_t^TH(s,\nu'(s),X^0(s),D V_{\nu'}(s,X^0(s)))ds-\sigma\int_t^T DV_{\nu'}(s,X^0(s))dB(s)\\\notag
&\quad+\int_t^T y(s,\nu'(s))dY^0(s)+\int_t^T r(s,\nu'(s))dR^0(s).
\end{align}
Let,
\begin{align}
\Delta V^n(t)&= V_{\nu^n}(t,X^0(t))-V_{\nu}(t,X^0(t)),\; \Delta DV^n(t)= DV_{\nu^n}(t,X^0(t))-DV_{\nu}(t,X^0(t)) \notag\\
\Delta g^n(T)&= \Delta V^n(T) = g(\nu^n(T),X^0(T))-g(\nu(T),X^0(T)),\notag\\
\Delta y^n(t)&= y(t,\nu^n(t))-y(t,\nu(t)),\;
\Delta r^n(t) = r(t,\nu^n(t))-r(t,\nu(t)), \notag\\
\Delta \psi^n(t)&= \sup_{\al\in U}|\psi(t,\nu^n(t),X^0(t),\al)-\psi(t,\nu(t),X^0(t),\al)|, \quad \psi \in \{f, b\}.
\label{eq_41}
\end{align}
Then, \eqref{eq_40} and \eqref{eq_41} imply,
\begin{align}\notag
&\Delta V^n(t)+ \sigma\int_t^T\Delta D V^n(s)dB(s)\\\notag
&\quad=\Delta g^n(T) + \int_t^T\Delta y^n(s)dY^0(s)+\int_t^T\Delta r^n(s)dR^0(s)\\\notag
&\qquad+\int_t^T[H(s,\nu^n(s),X^0(s),DV_{\nu^n}(s,X^0(s)))-H(s,\nu^n(s),X^0(s),DV_{\nu}(s,X^0(s)))]ds\\\notag
&\qquad+\int_t^T[H(s,\nu^n(s),X^0(s),DV_{\nu}(s,X^0(s)))-H(s,\nu(s),X^0(s),DV_{\nu}(s,X^0(s)))]ds.
\end{align}
Squaring both sides and then taking expectations gives
 \begin{align}\label{new200}
&\E[(\Delta V^n(t))^2]+ \sigma^2\E\left[\int_t^T(\Delta D V^n(s))^2ds\right]\\\notag
&\quad\le 2(T-t)\E\Big[\int_t^T[H(s,\nu^n(s),X^0(s),DV_{\nu^n}(s,X^0(s)))-H(s,\nu^n(s),X^0(s),DV_{\nu}(s,X^0(s)))]^2ds\Big]\\\notag
&\qquad +2C^n(t)\\\notag
&\quad\le 2c_B^2(T-t)\E\Big[\int_t^T(\Delta DV^n(s))^2ds\Big]+ 2C^n(t), 
\end{align}
where 
\begin{align}\notag
C^n(t)&=
\E\Big[ \Delta g^n(T) + \int_t^T\Delta y^n(s)dY^0(s)+\int_t^T\Delta r^n(s)dR^0(s)\\\notag
&\qquad\quad+\int_t^T[H(s,\nu^n(s),X^0(s),DV_{\nu}(s,X^0(s)))-H(s,\nu(s),X^0(s),DV_{\nu}(s,X^0(s)))]ds\Big]^2.
\end{align}
Letting $\delta=\sigma^2/(4c_B^2)$, we get from \eqref{new200} that for every $t\in[T-\delta,T]$,
\begin{align}\label{new201}
\E[(\Delta V^n(t))^2]+\frac{\sigma^2}{2}\E\left[\int_t^T(\Delta DV^n(s))^2ds\right]\le 4C^n(t). 
\end{align}
We now show that $\limsup_{n\to \infty} \sup_{0\le t \le T} C^n(t)=0$. 
Note that there exists $C_1\in (0, \infty)$ such that the following inequality holds for all $n \in \N$. 
\begin{align}\notag
\frac{1}{C_1} \sup_{0\le t \le T} C^n(t)&\le\E\left[(\Delta g^n(T))^2\right]+\E\left[(Y^0(T))^2\right]\sup_{0\le t\le T}(\Delta y^n(t))^2+\E\left[(R^0(T))^2\right]\sup_{0\le t\le T}(\Delta r^n(t))^2\\
&\quad+\E\left[
\int_0^T(\Delta f^n(s))^2ds\right]+\E\left[\int_0^T(\Delta b^n(s))^2ds \int_0^T (DV_\nu(s,X^0(s)))^2ds\right].\label{eq:eq1259}
\end{align}
Using the properties of the Skorohod map (Lemma \ref{lem_Skorohod}) it follows that 
$$ \E\left[(Y^0(T))^2 + (R^0(T))^2\right]<\iy .$$
The convergence of the right side of \eqref{eq:eq1259} to $0$ is now immediate from Assumption \ref{assumption1}(a), the boundedness of $DV_{\nu}$, and \eqref{eq_21}.
Thus from \eqref{new201} we have that
\begin{align}\label{new89}
\limn\E\left[\int_{T-\delta}^T|D V_{\nu^n}(t,X^0(t)) - D V_{\nu}(t,X^0(t))|dt\right]=0,\; \limn\E[(\Delta V^n(T-\delta))^2] =0.
\end{align}
Using the second convergence in \eqref{new89} and repeating the above argument for $t \in [T-\delta, T]$ to the interval $[T-2\delta, T-\delta]$, we see that
\eqref{new89} holds with $T$ replaced with $T-\delta$. Proceeding recursively in this manner we have \eqref{new79}. Hence we have shown that
\begin{equation}
	\label{eq:eq215}
	\al_{\nu^n}^0\to\al_\nu^0 \mbox{ in } \lambda_T\times\cal\PP \mbox{ - measure}.
\end{equation}
Using the above property we will now argue that $\Phi(\nu^n) \to \Phi(\nu)$ as $n\to \infty$, completing the proof of continuity of $\Phi$.
Let for $\nu' \in \clm_0$, 
$$\gamma_{\nu'}(t,x) = \hat \alpha(t, \nu'(t), x, DV_{\nu'}(t,x)), \;\; \hat b_{\nu'}(t,x) = b(t, \nu'(t), x, \gamma_{\nu'}(t,x)),\; (t,x) \in [0,T]\times [0,L]$$
and let $\PP^{\nu'}$ be a probability measure on $(\Omega, \clf)$ defined as
\begin{align}\label{up2}
d\PP^{\nu'} = \exp \left[ \frac{1}{\sigma} \int_0^T \hat b_{\nu'}(t, X^0(t)) dB(t) - \frac{1}{2\sigma^2} \int_0^T \hat b^2_{\nu'}(t, X^0(t)) dt\right] d\PP.
\end{align}
By Girsanov's theorem $\PP^{\nu'} \circ (X^0)^{-1} = \Phi(\nu')$. Thus to show $\Phi(\nu^n) \to \Phi(\nu)$ it suffices to argue that $\PP^{\nu^n} \to \PP^{\nu}$. We will in fact show that
$R(\PP^{\nu}\|\PP^{\nu^n}) \to 0$ as $n\to \infty$, where 
$$
R(\PP^{\nu}\|\PP^{\nu^n}) = \E^{\nu} \left(\log \frac{d\PP^{\nu}}{d\PP^{\nu^n}}\right) = \E \left( \frac{d\PP^{\nu}}{d\PP}\log \frac{d\PP^{\nu}}{d\PP^{\nu^n}}\right)$$
is the relative entropy of $\PP^{\nu}$ with respect to $\PP^{\nu^n}$, which due to Pinsker's inequality (see \cite[Page 132]{Tsybakov2009}) gives the convergence of $\PP^{\nu^n}$ to $\PP^{\nu}$.

Let
$$
\Delta \hat b_n(t,x) = \hat b_{\nu}(t,x) - \hat b_{\nu^n}(t,x), \;\; \Delta \hat b_n^2(t,x) = \hat b_{\nu}^2(t,x) - \hat b_{\nu^n}^2(t,x), \; (t,x) \in [0,T]\times [0,L].
$$
With this notation,
$$
\log \frac{d\PP^{\nu}}{d\PP^{\nu^n}} = \frac{1}{\sigma} \int_0^T \Delta \hat b_{n}(t, X^0(t)) dB(t) - \frac{1}{2\sigma^2} \int_0^T \Delta \hat b^2_{n}(t, X^0(t)) dt .$$
Also, noting that since $b$ is bounded $\E (\frac{d\PP^{\nu}}{d\PP})^2 \doteq \kappa < \infty$, we have from Cauchy-Schwarz inequality
\begin{equation}\label{eq:eq317}
R(\PP^{\nu}\|\PP^{\nu^n}) \le \sqrt{\kappa}\left(\E \left[\frac{1}{\sigma} \int_0^T \Delta \hat b_{n}(t, X^0(t)) dB(t) - \frac{1}{2\sigma^2} \int_0^T \Delta \hat b^2_{n}(t, X^0(t)) dt\right]^2\right)^{1/2} .\end{equation}
Next note that
\begin{align*}
	&\E \left[\frac{1}{\sigma} \int_0^T \Delta \hat b_{n}(t, X^0(t)) dB(t)\right]^2 \\
	&\quad = \frac{1}{\sigma^2} \int_0^T \E[\Delta \hat b_{n}(t, X^0(t))]^2 dt\\
	&\quad = \frac{1}{\sigma^2} \int_0^T \E \left[b(t, \nu(t), X^0(t), \alpha^0_{\nu}(t, X^0(t))) - b(t, \nu^n(t), X^0(t), \alpha^0_{\nu^n}(t, X^0(t))) \right]^2 dt.
\end{align*}
The last term converges to $0$ from the boundedness and continuity of $b$, \eqref{eq_21} and \eqref{eq:eq215}. Similarly
$$
\E \left[\frac{1}{2\sigma^2} \int_0^T \Delta \hat b^2_{n}(t, X^0(t)) dt\right]^2 \to 0$$
as $n \to \infty$. Using the above two observations in \eqref{eq:eq317} we have $R(\PP^{\nu}\|\PP^{\nu^n}) \to 0$ and thus,
as argued earlier, the proof of continuity of $\Phi$ is complete. Thus we have shown that $\Phi$ is a continuous map on $\calM_0$, which is a non-empty convex subset of the normed space $\calB$, into a compact set $\calK\subset\calM_0$. Thus by the fixed point theorem in Lemma \ref{lem_Schauder}, $\Phi$ has a fixed point.\\

The first results on unique solvability of a MFG go back to \cite{Lasry2007}. Since then uniqueness has been argued in various settings (see e.g.~\cite{Cardaliaguet2013,Gomes2013,Carmona2015}). The proof given below uses arguments similar to those in \cite[Section 7.3]{Carmona2015}, however for the sake of completeness we give the details. 
Consider as before a system $\Xi =(\Omega,\calF,\{\calF_t\},\PP, B)$ and let $Z^0 = (X^0, Y^0, R^0)$ be given through \eqref{eq:eq433}.
Let $\nu_1,\nu_2\in\calM_0$. For $i=1,2$, let $\al^0_{\nu_i}$ and $\PP^{\nu_i}$ be given by \eqref{up1} and \eqref{up2}, respectively, with $\nu_i$ replacing $\nu'$. 
Applying It\^{o}'s lemma to $V_{\nu_1}(t,X^0(t))$ and recalling that $V_{\nu_1}$ solves \eqref{HJB1}--\eqref{HJB2} with $\nu$ replaced with $\nu_1$, we get
\begin{align*}
	&V_{\nu_1}(T,X^0(T))-V_{\nu_1}(0,x) \\
	&\quad = - \int_0^T H(t, \nu_1(t), X^0(t), DV_{\nu_1}(t,X^0(t))) dt
	+ \sigma \int_0^T DV_{\nu_1}(t,X^0(t)) dB(t) + \zeta_T\\
	&= - \int_0^T h(t, \nu_1(t), X^0(t), \alpha^0_{\nu_1}(t), DV_{\nu_1}(t,X^0(t))) dt
	+ \sigma \int_0^T DV_{\nu_1}(t,X^0(t)) dB(t) + \zeta_T,\\
\end{align*}
where recalling the form of $y$ and $r$ from Assumption \ref{assumptionU} 
\begin{align*}
	\zeta_T &= \int_{[0,T]} DV_{\nu_1}(t,0) dY^0(t) + \int_{[0,T]} DV_{\nu_1}(t,L) dR^0(t)\\
&= \int_{[0,T]} y(t) dY^0(t) + \int_{[0,T]} r(t) dR^0(t).
\end{align*}
Observing that for $t \in [0,T]$
$$
h(t, \nu_1(t), X^0(t),\alpha^0_{\nu_1}(t), DV_{\nu_1}(t,X^0(t)))
= h(t, \nu_1(t), X^0(t),\alpha^0_{\nu_1}(t), DV_{\nu_2}(t,X^0(t))) + \clr_t,$$
where
$$\clr_t = [DV_{\nu_1}(t,X^0(t)) - DV_{\nu_2}(t,X^0(t))] \hat b_{\nu_1}(t, X^0(t)),$$
we have
\begin{align}\label{asaf1}
	V_{\nu_1}(0,x) &= g(\nu_1(T),X^0(T)) + \int_0^T h(t, \nu_1(t), X^0(t), \alpha^0_{\nu_1}(t), DV_{\nu_2}(t,X^0(t))) dt\\\notag
	&\quad + \int_0^T \clr_t dt - \sigma \int_0^T DV_{\nu_1}(t,X^0(t)) dB(t) - \zeta_T .
\end{align}
Similarly, applying It\^{o}'s lemma to $V_{\nu_2}(t,X^0(t))$
\begin{align}\label{asaf2}
	V_{\nu_2}(0,x) &= g(\nu_2(T),X^0(T)) + \int_0^T H(t, \nu_2(t), X^0(t), DV_{\nu_2}(t,X^0(t))) dt\\\notag
	&\quad - \sigma \int_0^T DV_{\nu_2}(t,X^0(t)) dB(t) - \zeta_T.
\end{align}
Substracting \eqref{asaf2} from \eqref{asaf1}
\begin{align*}
	&V_{\nu_1}(0,x) - 	V_{\nu_2}(0,x) \\
	&\quad = g(\nu_1(T),X^0(T)) - g(\nu_2(T),X^0(T))
	- \sigma \int_0^T [DV_{\nu_1}(t,X^0(t))- DV_{\nu_2}(t,X^0(t))] dB^{\nu_1}(t)\\
	&\quad + \int_0^T \left[h(t, \nu_1(t), X^0(t), \alpha^0_{\nu_1}(t), DV_{\nu_2}(t,X^0(t)))
	- H(t, \nu_2(t), X^0(t), DV_{\nu_2}(t,X^0(t)))\right] dt,
\end{align*}
where for $i=1,2$, $B^{\nu_i}(t)\doteq B(t)-\frac{1}{\sigma}\int_0^t\hat b_{\nu_i}(s,X^0(s))ds$, $t\in [0,T]$. Since under $\PP^{\nu_i}$, $B^{\nu_i}$ is a standard Brownian motion, taking expectation under the measure $\PP^{\nu_1}$
%
%
%
%
\begin{align}\notag
&V_{\nu_1}(0,x)-V_{\nu_2}(0,x)\\
&\quad=\E^{\nu_1}\Big[\int_0^T\left[
h(t,\nu_1(t),X^0(t),\al^0_{\nu_1}(t),DV_{\nu_2}(t,X^0(t)))-H(t,\nu_2(t),X^0(t),DV_{\nu_2}(t,X^0(t)))\right]dt\notag\\
&\qquad\qquad\qquad+g(\nu_1(T),X^0(T))-g(\nu_2(T),X^0(T))
\Big].\label{up3}
\end{align}
A similar calculation shows
\begin{align}\notag
&V_{\nu_1}(0,x)-V_{\nu_2}(0,x)\\
&\quad=\E^{\nu_2}\Big[\int_0^T\left[
H(t,\nu_1(t),X^0(t),DV_{\nu_1}(t,X^0(t)))-h(t,\nu_2(t),X^0(t),\al^0_{\nu_2}(t),DV_{\nu_1}(t,X^0(t)))\right]dt \notag\\
&\qquad\qquad\qquad+g(\nu_1(T),X^0(T))-g(\nu_2(T),X^0(T))
\Big].\label{up3b}
\end{align}
By the definition of the $H$ and the form of $f$ in Assumption \ref{assumptionU} we get,
\begin{align}\notag
&h(s,\nu_1(s),X^0(s),\al^0_{\nu_1}(s),DV_{\nu_2}(s,X^0(s)))\\\notag
&\quad\ge H(s,\nu_1(s),X^0(s),DV_{\nu_2}(s,X^0(s)))\\\notag
&\quad=f_0(s,\nu_1(s),X^0(s))+f_1(s,X^0(s),\al^0_{\nu_2}(s))+DV_{\nu_2}(s,X^0(s))\hat b_{\nu_2}(s)
\end{align}
where the last equality uses the observation that since
$\hat \alpha(s, \eta,x,p)$ does not depend on $\eta$, 
$\hat \alpha(s, \nu_1(s), X^0(s), DV_{\nu_2}(s,X^0(s))) = \alpha^0_{\nu_2}(s)$.
Therefore for all $s\in [0,T]$, 
\begin{align}
&h(s,\nu_1(s),X^0(s),\al^0_{\nu_1}(s),DV_{\nu_2}(s,X^0(s)))-H(s,\nu_2(s),X^0(s),DV_{\nu_2}(s,X^0(s))) \notag\\
&\quad\ge f_0(s,\nu_1(s),X^0(s))-f_0(s,\nu_2(s),X^0(s)), \label{up5}
\end{align}
Similarly for all $s\in [0,T]$,
\begin{align}
&H(s,\nu_1(s),X^0(s),DV_{\nu_1}(s,X^0(s)))-h(s,\nu_2(s),X^0(s),\al^0_{\nu_2}(s),DV_{\nu_1}(s,X^0(s)))\notag\\
&\quad\le f_0(s,\nu_1(s),X^0(s))-f_0(s,\nu_2(s),X^0(s)). \label{up6}
\end{align}
Applying the last two inequalities to \eqref{up3} and \eqref{up3b}, we get
\begin{align}\label{newnew1}
0\le [\E^{\nu_2}-\E^{\nu_1}]\Big[&g(\nu_1(T),X^0(T))-g(\nu_2(T),X^0(T))\\\notag
&\quad+\int_0^T[f_0(s,\nu_1(s),X^0(s))-f_0(s,\nu_2(s),X^0(s))]ds\Big].
\end{align}
Until now $\nu_1$ and $\nu_2$ were arbitrary measures in $\clm_0$. Suppose now that $\nu_i$, $i=1,2$, are fixed points
of $\Phi$. Then, for $i=1,2$, $\PP^{\nu_i} \circ (X^0)^{-1} = \nu_i$ and so for all $s \in [0,T]$,
$\PP^{\nu_i} \circ (X^0(s))^{-1} = \nu_i(s)$.
In this case, using the inequalities in Assumption \ref{assumptionU} we get that the inequality \eqref{newnew1} can be replaced with equality. We claim that $\al^0_{\nu_1}=\al^0_{\nu_2}$, $\lambda_T^0\times\PP$-a.s. Indeed, suppose that there is a set with positive $\lambda_T^0\times\PP$-measure on which $\al^0_{\nu_1}\neq\al^0_{\nu_2}$. Then on this set \eqref{up5} and \eqref{up6} will hold with strict inequalities by Assumption \ref{assumption1} (b). Since the measures $\PP$, $\PP^{\nu_1}$, and $\PP^{\nu_2}$ are equivalent, this will say that \eqref{newnew1} holds with a strict inequality as well, which contradicts the equality that was established above. This proves the claim. Since $b(t,\eta,x,u)$ does not depend on $\eta$ we conclude from the equality of $\al^0_{\nu_1}$ and $\al^0_{\nu_2}$ that $\hat b_{\nu_1}(t, X^0(t,\omega)) = \hat b_{\nu_2}(t, X^0(t,\omega))$, $\lambda_T^0\times\PP$-a.s. and thus by
\eqref{up2}, $\PP^{\nu_1}=\PP^{\nu_2}$.
Combining this with the fact that $\nu_i$ are fixed points of $\Phi$, we now have
 $\nu_1=\PP^{\nu_1}\circ (X^0)^{-1}=\PP^{\nu_2}\circ (X^0)^{-1}=\nu_2$. \hfill$\Box$

\section{Asymptotic Nash equilibrium}\label{sec4}
\beginsec
The main result of this section is Theorem \ref{thm1}. The main idea in the proof is to use a solution $\bar \nu$ to the MFG (which from Theorem \ref{thm_fixed} exists under Assumption \ref{assumption1}) and the associated feedback control given by Proposition \ref{lem2} in order to construct an admissible control
$\tilde \alpha^n = \{\tilde \alpha^{n,i}\}_{i=1}^n$ for the $n$-player game. Specifically, the control
will be given in a feedback form through the following relation
\begin{align}
	\tilde \alpha^{n,i}(t) \doteq \hat \al(t,\bar\nu(t),\tilde Q^{n,i}(t), D V_{\bar \nu}(t,\tilde Q^{n,i}(t))),
\label{eq:eq424}	
\end{align}
where $\tilde Q^{n,i}$ is the corresponding scaled queue length under the feedback control. Note that the only information each of the players uses is its own state. Therefore, the problem is decentralized in the sense that players do not need to observe each others' states. 
\color{black} 
Our main condition, in addition to Assumption \ref{assumption1}, for $\{\tilde \alpha^n\}$ to be an asymptotic Nash equilibrium is the following. It in particular says that the drift function does not depend on the mean-field term.
Fix $\eta_0 \in \clp([0,L])$.
\begin{assumption}\label{assumption4}
$\,$
\begin{enumerate}
\item[(a)]	
For every $t\in[0,T]$, $\eta\in\calP([0,L])$, $x\in[0,L]$, and $\al\in U$, one has,
\begin{align}\notag
b(t,\eta,x,\al)&=b(t,\eta_0,x,\al);
\end{align}
\item[(b)] There exists $x\in[0,L]$ such that for every $i\in\N$,
\begin{align}
\limn\tQni(0)=x.
\end{align}
\end{enumerate}
\end{assumption}
As before, with an abuse of notation, we will write $b(t,\eta,x,\al)$ as $b(t,x,\al)$ when Assumption \ref{assumption4} holds.
As discussed in Remark \ref{rem_uniqueness}, part (a) of the assumption means that the empirical measure affects the drift only through the control, which in turn is affected by the empirical measure through the running cost. 

Recall the probability space $(\Omega', \clf', \PP')$ from Section \ref{sec2a}.
Let for $n \in \N$, $t \in [0,T]$ and $i = 1, \cdots ,n$, $\beta^n(t): \Omega' \to U$ be such that
$$\tilde \alpha^n_{-i} =\{ \tilde \alpha^{n,1}, \ldots , \tilde \alpha^{n, i-1}, \beta^n, \tilde \alpha^{n, i+1}, \ldots , \tilde \alpha^{n,n}\}$$
is an admissible control (i.e.~$\tilde \alpha^n_{-i} \in \clu^n$).
The following is the main result of the section.

\begin{theorem}\label{thm1}
Suppose Assumptions \ref{assumption1} and \ref{assumption4} hold.
Let $\bar \nu$, $\tilde \alpha^n$ and $\tilde \alpha^n_{-i}$, $i=1, \ldots, n$, $n \in \N$ be as introduced above.
Then
\begin{align}\label{20aa}
\limsup_{n\to \infty} J^{n,1}(\tQn(0);\tilde \al^{n})= V_{\bnu}(0,x)\le\liminf_{n\to \infty} J^{n,1}(\tQn(0);
\tilde \al^{n}_{-i}).
\end{align}
\end{theorem}
Theorem \ref{thm1} in particular says that for every $\eps>0$, there is $n\in\N$ sufficiently large such that $\tilde \alpha^n$ forms an $\eps$-Nash equilibrium in the $n$-player game.

The proof is given in the next three sections. First in Section \ref{4z} (Proposition \ref{prop:tightemp}) we will prove
the convergence of empirical measures of the scaled queue length processes under controls 
$\tilde \al^{n}$ and $\tilde \al^{n}_{-i}$ to $\bar \nu$. Next, in Section \ref{sec4a} (Proposition \ref{prop2}) we will prove the first equality in \eqref{20aa} and finally Proposition \ref{prop3} in Section \ref{sec4b} will prove the 
inequality in \eqref{20aa}.


\subsection{Convergence of empirical measures}\label{4z}
 Let for $i\in \N$, $\tilde\alpha^n_{-i} \in \clu^n$ be as defined below Assumption \ref{assumption4}.
Let
$$\tilde \nu^n_{-i} = \frac{1}{n} \sum_{j=1}^n \delta_{\tilde Q^{n,j}},$$
where $\tilde Q^{n,j}$ is the controlled queue length process defined by \eqref{7} with $\alpha^n$ replaced with
$\tilde \alpha^n_{-i}$. The following result gives the convergence of $\tilde \nu^n_{-i}$ to $\bar \nu$.
\begin{proposition}\label{prop:tightemp}
	Suppose Assumptions \ref{assumption1} and \ref{assumption4} hold. Then for every $i \in \N$,
	$\tilde \nu^n_{-i}$ converges in probability, in $\clp(\cld([0,T]:[0,L]))$, to $\bar \nu$ as $n \to \infty$.
\end{proposition}
{\bf Proof:} 
Without loss of generality we assume that $i=1$. Recall that for $j \in \N$, $\tilde Q^{n,j}$ is defined by \eqref{7}
with $\alpha^n$ replaced with $\tilde \alpha^n_{-1}$. Define 
$$\tilde\zeta^{n,i}(t) = \tilde Q^{n,i}(0) + \tilde A^{n,i}(t) - \tilde D^{n,i}(t) + \int_0^t\tbni(s)ds, \; t \in [0,T],$$
where $\tilde A^{n,i}, \tilde D^{n,i}$ are as in \eqref{eq:eq928}. Define for $i = 1, \ldots, n$,
$$\tilde G^{n,i} \doteq (\zeta^{n,i}, \tilde Q^{n,i}, \tilde Y^{n,i}, \tilde R^{n,i})$$
and let
$$\Xi^n \doteq \frac{1}{n-1} \sum_{i=2}^n \delta_{\tilde G^{n,i}}.$$
Note that since by assumption $\{\tilde Q^{n,i}(0)\}_{i=1}^n$ are exchangeable and the controls $\tilde \alpha^{n,i}$
are given in terms of the same feedback function $\hat \alpha$ for each $i= 2, \ldots , n$, the processes
$\{\tilde G^{n,i}\}_{i=2}^n$ are exchangeable. Defining $\tilde G^{n,i}$ to be the zero process for $i > n$
we can regard, $\tilde G^n \doteq \{\tilde G^{n,i}\}_{i=2}^{\infty}$ as a random variable with values in
$\cld([0,T]: (\R^4)^{\otimes\infty})$. We now argue the tightness of the sequence $\{\tilde G^n\}$.
It suffices to show for each $i$, the tightness of $\{\tilde G^{n,i}\}_{n\in\N}$ in $\cld([0,T]: \R^4)$.
Since 
\begin{equation}
	\sup_{n,i,\omega} \sup_{t\in [0,T]} \frac{\la^{n,i}(t,\omega)+ \mu^{n,i}(t,\omega)}{e_n} \doteq C_0 < \infty, \label{eq:eq858}
\end{equation}
the following two conditions are satisfied with $X^n$ equal to $\langle\tAni, \tAnj \rangle$ and $\langle\tDni, \tDnj \rangle$ for all $i,j$.
\begin{description}
\item{[A]} For each $\eps>0,\eta>0$ there exists a $\del>0$ and $n_0\in \N$ with the property that for every family of stopping times $\{\tau_n\}_{n\in\N}$ ($\tau_n$ being an $\clf_t$-stopping time on $(\Omega', \clf', \PP')$) with $\tau_n\leq T-\delta$,
\begin{align*}
\sup_{n\geq n_0}\sup_{\theta\leq \del}\PP\{\|X^n(\tau_n)-X^n(\tau_n+\theta)\|\geq \eta\}\leq \eps.
\end{align*}
\item{[T$_1$]} For every $t$ in some dense subset of $[0,T]$, $\{X^n(t)\}_{n\in\N}$ is a tight sequence of $\R$ valued random variables.
\end{description}
Then by Rebolledo's theorem (see \cite[Theorem 2.3.2]{joffe1986weak})
$\{\tilde A^{n,i}\}_{n\ge 1}$ and $\{\tilde D^{n,i}\}_{n\ge 1}$ are tight in $\cld([0,T]:\R)$ for each $i$.
Also since the jumps of $\tilde A^{n,i}$ and $\tilde D^{n,i}$ are of size $1/\sqrt{e_n}$, these processes are
$\clc$-tight (namely all weak limit points are continuous a.s.). From boundedness of $b$ we see that
$\{\int_0^\cdot\tbni(s)ds\}_{n\in \N}$ is tight in $\clc([0,T]:\R)$.
Combining this with Assumption \ref{assumption4}(b) we see that $\{\zeta^{n,i}\}_{n\in \N}$ is $\clc$-tight
in $\cld([0,T]:\R)$. Using now the continuity of the Skorohod map (Lemma \ref{lem_Skorohod}) we have the desired tightness of
$\{\tilde G^{n,i}\}_{n\in\N}$. 

Suppose now that, along some subsequence, $\tilde G^n$ converges to $\tilde G \doteq (\zeta^{i}, \tilde Q^{i}, \tilde Y^{i}, \tilde R^{i})$, in distribution, in
$\cld([0,T]: (\R^4)^{\otimes\infty})$. Then $\tilde G \in \clc([0,T]: (\R^4)^{\otimes\infty})$ a.s. and from \cite[Lemma 4.2]{Kotelenez2010} and the exchangeability of $\{\tilde G^{n,i}\}_{i=2}^n$ it follows that 
$\{\tilde G^{i}\}_{i=2}^{\infty}$ is exchangeable and(along the subsequence),
\begin{equation}
	\label{eq:eq453}
	(\tilde G^n, \Xi^n) \Rightarrow (\tilde G, \Xi)
\end{equation}
in $\cld([0,T]: (\R^4)^{\otimes\infty})\times \clp(\cld([0,T]:\R^4))$ where 
$\Xi \doteq \lim_{m\to \infty} \frac{1}{m-1}\sum_{i=2}^m \delta_{\tilde G^i}$.

We will now characterize the distribution of $\{\tilde Q^i\}$. From tightness of $\{\tilde Y^{n,i}\}_{n\in N}$
and $\{\tilde R^{n,i}\}_{n\in N}$ argued above and \eqref{8} it follows that
$$ \frac{1}{e_n}\int_0^t 1_{\{\tQni(s)=0\}}\muni(s)ds \to 0,\;\mbox{ and }
\frac{1}{e_n}\int_0^t 1_{\{\tQni(s)=L\}}\lani(s)ds \to 0,$$
uniformly on $[0,T]$, in probability.
Also, from \eqref{2} it follows that
$$\sup_{0\le t \le T} \left[ \left| \frac{\la^{n,i}(t)}{e_n} - \hat \lambda\right| +
\left| \frac{\mu^{n,i}(t)}{e_n} - \hat \mu\right|\right] \to 0 \mbox{ a.s. }$$
as $n\to \infty$.
Thus from \eqref{eq: eq415} (and the relation $\hat \lambda = \hat \mu$), for all $i,j$,
$$
\langle\tAni, \tAnj \rangle(t) \to \delta_{ij}\hat \lambda,\; \langle \tDni, \tDnj \rangle(t) \to
\delta_{ij}\hat \lambda,\; \langle\tAni, \tDnj \rangle(t) \to 0$$
in probability, uniformly on $[0,T]$, as $n\to \infty$.
By standard martingale techniques it now follows that
$$\{\tAni - \tDni\}_{i\ge 1} \Rightarrow \{\sigma B^i\}_{i\ge 1},$$
in $D([0,T]:\R^{\infty})$, where $\{B^i\}$ are mutually independent standard Brownian motions.
Also, since for $i\ge 2$
$$
(\tQni,\tYni,\tRni)(t)=\Gamma\left(\tQni(0)+\int_0^\cdot\tbni(s)ds+\tAni(\cdot)-\tDni(\cdot) + o(1)\right)(t), \; t \in [0,T],$$
where
$$\tbni(t) = b(t,\tQni(t), \hat \al(t,\bar\nu(t),\tilde Q^{n,i}(t), D V_{\bar \nu}(t,\tilde Q^{n,i}(t)))),$$
we have from the continuity of $b$ (Assumption \ref{assumption1}), $\hat \alpha$, and $D V_{\bar \nu}$, for $i\ge 2$,
\begin{equation}\label{eq:eq637}
(\tilde Q^i,\tilde Y^i,\tilde R^i)(t)=\Gamma\left(x+\int_0^\cdot
b(t,\tilde Q^i(t), \hat \al(t,\bar\nu(t),\tilde Q^{i}(t), D V_{\bar \nu}(t,\tilde Q^{i}(t))))dt +\sigma B^i(\cdot) \right)(t).
\end{equation}
Once again using standard martingale arguments it follows that for $0\le s \le t \le T$,
$B^i(t)-B^i(s)$ is independent of $\sigma\{(\tilde Q^i(u), \tilde R^i(u), \tilde Y^i(u), B^i(u)): u \le s\}$.
From weak uniqueness property noted in Section \ref{sec3b1} and the fact that $\bar \nu$ is a fixed point of $\Phi$ we now have that $\tilde Q^i$ has distribution $\bar \nu$ for $i=2,3, \cdots$. 
 Using the fact that $\{B^i\}$ are mutually independent, a simple argument based on Girsanov's theorem shows that $\{\tilde Q^i\}$ are mutually independent as well. This characterize the distribution
of $\{\tilde Q^i\}_{i\ge 2}$ as $\bar \nu^{\otimes \infty}$. We now have from \eqref{eq:eq453}, the definition of $\Xi$, and the law of large numbers that
$$\lim_{n\to \infty }\tilde \nu^n_{-i} = \lim_{n\to \infty}\frac{1}{n} \sum_{i=1}^n \delta_{\tilde Q^{n,i}} 
= \lim_{n\to \infty}\frac{1}{n-1} \sum_{i=2}^n \delta_{\tilde Q^{n,i}} = \lim_{m\to \infty} \frac{1}{m-1}\sum_{i=2}^m \delta_{\tilde Q^i} = \bar \nu.$$
The result follows. \hfill \qed
\begin{remark}
	\label{rem:rem640}
	The above proof also shows that if $\tilde \alpha^n_{-1} = \tilde \alpha^n$, namely $\beta^n = \tilde \alpha^{n,1}$, then \eqref{eq:eq637}
	holds for all $i\ge 1$ and the law of $\{\tilde Q^i\}_{i\ge 1}$ is $\bar \nu^{\otimes \infty}$.
\end{remark}

\subsection{Same strategy for all players}\label{sec4a}
In this section we prove the 
equality in \eqref{20aa}. 

\begin{proposition}\label{prop2}
Suppose Assumptions \ref{assumption1} and \ref{assumption4} are satisfied. Let $\tilde \alpha^n = \{\tilde \alpha^{n,i}\}_{i=1}^{\infty}$ be as in \eqref{eq:eq424}. Then for all $i\ge 1$	
\begin{align}\label{20}
\limn J^{n,i}(\tQn(0);\tilde \al^{n})= V_{\bnu}(0,x).
\end{align}
\end{proposition}
\noi
{\bf Proof:} 
Without loss of generality we assume $i=1$. From the proof of Proposition \ref{prop:tightemp} (see Remark \ref{rem:rem640})
$$(\tilde Q^{n,1}, \tilde Y^{n,1}, \tilde R^{n,1}) \Rightarrow (\tilde Q^{1}, \tilde Y^{1}, \tilde R^{1})$$
where the processes on the right side are given through \eqref{eq:eq637} with $i=1$. 
Let
$$\hat \alpha(t, \bar \nu(t), \tilde Q^1(t) , D V_{\bar \nu}(t,\tilde Q^{1}(t))) \doteq \gamma(t, \tilde Q^1(t)), \; t \in [0,T].$$
Recall that $f$ and $g$ are bounded continuous functions and from Proposition \ref{prop:tightemp} we have that, for every $t \in [0,T]$,
$$
(\tilde \nu^n(t), \tilde Q^{n,1}(t), \alpha^{n,1}(t)) \Rightarrow (\bar \nu(t), \tilde Q^{1}(t), \gamma(t, \tilde Q^1(t))).$$
This shows that
\begin{align}
&\E\left[\int_0^T f(t,\tnun(t), \tilde Q^{n,1}(t), \alpha^{n,1}(t))dt + g(\tnun(T),\tilde Q^{n,1}(T))\right]\nonumber\\
	&\quad \to \E\left[\int_0^T f(t,\bar \nu(t), \tilde Q^1(t), \gamma(t, \tilde Q^1(t)))dt + g(\bar \nu(T),\tilde Q^1(T))\right].\label{eq:eq656}
\end{align}
Also by continuity of $y$ and $r$
$$
\left( y(\cdot, \tilde \nu^n(\cdot)), r(\cdot, \tilde \nu^n(\cdot)), \tilde Y^{n,1}(\cdot), \tilde R^{n,1}(\cdot)\right)
\Rightarrow \left( y(\cdot, \bar \nu(\cdot)), r(\cdot, \bar \nu(\cdot)), \tilde Y^{1}(\cdot), \tilde R^{1}(\cdot)\right)$$
in $D([0,T]:\R^4)$.
It then follows (cf.~\cite[Lemma 2.4]{dai1996existence}) 
\begin{equation}
	\label{eq:eq853}
\left ( \int_0^Ty(t,\tnun(t))d\tilde Y^{n,1}(t), \int_0^Tr(t,\tnun(t))d\tilde R^{n,1}(t)\right)
\Rightarrow \left ( \int_0^Ty(t,\bar \nu(t))d\tilde Y^1(t), \int_0^Tr(t,\bar \nu(t))d\tilde R^1(t)\right)\end{equation}
as $n\to \infty$.
Also from Lemma \ref{lem_Skorohod},
\begin{equation}\label{eq:eq905}\E\left [ (\tilde Y^{n,1}(T))^2 + (\tilde R^{n,1}(T))^2 \right] \le c_S^2 \E\left [\sup_{0\le s \le T} |\tilde \zeta^{n,1}(s)|^2\right].\end{equation}
Next note that
$$\E \sup_{0\le t \le T} (\tilde A^{n,1}(t))^2 \le 4 \frac{1}{e_n} \int_0^T \la^{n,1}(s) ds \le 4C_0T,$$
where $C_0$ is as in \eqref{eq:eq858}.
Similarly, $\E \sup_{0\le t \le T} (\tilde D^{n,1}(t))^2 \le 4C_0T$.
Combining these estimates
$$\sup_{n}\E \left[\sup_{0\le t \le T}|\tilde \zeta^{n,1}(s)|^2\right] < \infty$$
which combined with \eqref{eq:eq853} and the boundedness of $y,r$ implies
$$\sup_{n}\E\left[\int_0^Ty(t,\tnun(t))d\tilde Y^{n,1}(t) + \int_0^Tr(t,\tnun(t))d\tilde R^{n,1}(t) \right]^2 < \infty.$$
Combining this with the weak convergence in \eqref{eq:eq853} we have
\begin{align}
	&\E\left[\int_0^Ty(t,\tnun(t))d\tilde Y^{n,1}(t) + \int_0^Tr(t,\tnun(t))d\tilde R^{n,1}(t) \right]\nonumber \\
	&\quad
	\to \E\left[\int_0^Ty(t,\bar \nu(t))d\tilde Y^1(t)+ \int_0^Tr(t,\bar \nu(t))d\tilde R^1(t)\right].\label{eq:eq908}
\end{align}
Combining \eqref{eq:eq656} and \eqref{eq:eq908} and recalling from the optimality of $\hat \alpha$ that
\begin{align}\label{new5}
V_{\bnu}(0,x)=\E\Big[&\int_0^T f(t,\bnu(t),\tilde Q^1(t),\hat\al(t,\bnu(t),\tilde Q^1(t), DV_{\bnu}(t,\tilde Q^1)))dt + g(\bnu(T),\tilde Q^1(T))\\\notag
&\quad+\int_0^Ty(t,\bnu(t))d\tilde Y^1(t)+\int_0^Tr(t,\bnu(t))d\tilde R^1(t)\Big],
\end{align}
we have the desired convergence $\limn J^{n,1}(\tQn(0);\tilde \al^{n})= V_{\bnu}(0,x)$. \hfill \qed

\subsection{Deviation of Player 1}\label{sec4b}
In this section we prove the inequality on the right side of \eqref{20aa}. 
\begin{proposition}\label{prop3}
Suppose Assumptions \ref{assumption1} and \ref{assumption4} hold.
Let $\bar \nu, \beta^n, \tilde \alpha^n_{-i}$ be as introduced at the beginning of Section \ref{sec4}. Then for each $i\ge 1$
\begin{align}\label{500}
\limninf J^{n,1}(\tQn(0);\tilde \al^{n}_{-i})\ge V_{\bnu}(0,x).
\end{align}
\end{proposition}
{\bf Proof:} 
As before, we assume without loss of generality that $i=1$. We will need to argue the tightness of the control sequence $\{\beta^n\}$ in an appropriate space.
For this it will be convenient to consider a relaxed control formulation. Consider the relaxation of the stochastic control problem in \eqref{13}--\eqref{15} where the control space $U$ is replaced by
$\clp(U)$, the drift function $b$ is replaced by the function $b_{\clr}:[0,T]\times [0,L]\times \clp(U) \to \R$
defined as
$$b_{\clr}(t,x, r) \doteq \int_U b(t, x, u) r(du), \; (t, x, r)\in [0,T]\times [0,L]\times \clp(U),$$
and the running cost function $f$ is replaced by $f_{\clr}: [0,T]\times \clp([0,L])\times [0,L]\times \clp(U) \to \R$, defined as
$$f_{\clr}(t,\eta, x, r) \doteq \int_U f(t, \eta, x, u) r(du), \; (t, \eta, x, r)\in [0,T]\times \clp([0,L])\times [0,L]\times \clp(U).$$
Also, the class of admissible controls $\cla(\Xi, t,x, \bar \nu)$ is replaced by $\cla_{\clr}(\Xi, t,x, \bar \nu)$ of pairs $(\alpha_{\clr}, Z)$ that are similar to pairs 
$(\alpha, Z)$ introduced above \eqref{13} except that $\alpha_{\clr}$ is $\clp(U)$ valued rather than $U$ valued and \eqref{13} holds with
$\bar b(u) = b(u, X(u), \alpha(u))$ replaced with $b_{\clr}(u, X(u), \alpha_{\clr}(u))$. The corresponding cost function
$J_{\bar \nu, \clr}$ is defined by \eqref{eq:eq234} with $f$ replaced by $f_{\clr}$. The value function in this relaxed formulation, denoted as $V_{\bar \nu, \clr}$,
is given by \eqref{15} with $\cla$ replaced by $\cla_{\clr}$. Define the function $h_{\clr}$ by \eqref{19}, replacing $(f,b)$ with $(f_{\clr},b_{\clr})$.
Then, from Assumption \ref{assumption1}(b),
$$H(t,\eta,x,p) = \inf_{u \in U} h(t,\eta, x, u,p) = \inf_{r \in \clp(U)} h_{\clr}(t,\eta, x, r,p).$$
This shows that $V_{\nu}$ and $V_{\nu, \clr}$ are both solutions of the PDE \eqref{HJB1} - \eqref{HJB2}. In view of the uniqueness result from Proposition \ref{lem2}, $V_{\nu} = V_{\nu, \clr}$.

Let $\beta^n_{\clr}(t) \doteq \delta_{\beta^n(t)}$, $t \in [0,T]$ and define $\bar \beta^n_{\clr} \in \clm(U \times [0,T])$ as
$$\bar \beta^n_{\clr} (du\, dt) \doteq \beta^n_{\clr}(t)(du) dt,$$
where $\clm(U \times [0,T])$ is the space of finite measures on $U \times [0,T]$ equipped with the topology of weak convergence.
Then we can rewrite
\begin{align}\label{eq:eq300}
J^{n,1}(\tilde Q^n(0), \tilde \alpha^n_{-1})\doteq &
\E\Big[\int_{U\times [0,T]} f(t,\tnun(t), \tilde Q^{n,1}(t), u) \bar \beta^n_{\clr}(du\,dt) + g(\tnun(T),\tilde Q^{n,1}(T))\\\notag
&\quad+\int_0^Ty(t,\tnun(t))d\tilde Y^{n,1}(t)+\int_0^Tr(t,\tnun(t))d\tilde R^{n,1}(t)\Big]
\end{align}
and
\begin{align}\label{eq:eq303}
(\tilde Q^{n,1},\tilde Y^{n,1},\tilde R^{n,1})=\Gamma\left(\tilde Q^{n,1}(0)+
\int_{U\times [0,\cdot]}b(s, \tilde Q^{n,1}(s), u) \bar \beta^n_{\clr}(du\,ds) +\tAni(\cdot)-\tDni(\cdot) + o(1)\right).
\end{align}
From Proposition \ref{prop:tightemp}, $\tnun \Rightarrow \bar \nu$. Also, the arguments of the same proposition show that
\begin{equation}
	\label{eq:eq426}\left\{\tilde Q^{n,1}(\cdot), \tilde A^{n,1}(\cdot) - \tilde D^{n,1}(\cdot), \int_{U\times [0,\cdot]} b(s, \tilde Q^{n,1}(s), u) \bar \beta^n_{\clr}(du\,ds) \right\}_{n\in \N}
\end{equation}
are $\clc$-tight in $\cld([0,T]: \R^3)$. Furthermore, since $U\times [0,T]$ is compact and $\bar \beta^n_{\clr}(U\times [0,T]) = T$, the sequence
$\{\bar \beta^n_{\clr}\}_{n\in \N}$ is tight in $\clm(U\times [0,T])$.
Suppose now that along a subsequence (labeled once more as $\{n\}$) the sequence in \eqref{eq:eq426} along with $\{\bar \beta^n_{\clr}\}$ converges in distribution to
$(\tilde Q^1, \sigma B^1, \vartheta, \bar \beta_{\clr}).$
Then from the Lipschitz property of $b$ (Assumption \ref{assumption1}(a)) we have, for $t\in [0,T]$,
\begin{align*}
\vartheta(t) &= \int_{U\times [0,t]} b(s, \tilde Q^1(s), u) \bar \beta_{\clr}(du ds)\\
&= \int_{0}^t b_{\clr}(s, \tilde Q^1(s), \beta_{\clr}(s)) ds,
\end{align*}
where $\beta_{\clr}(s)$ is obtained by disintegrating $\bar \beta$, i.e.~$\bar \beta_{\clr}(du\, ds) = \beta_{\clr}(s)(du) ds$.
Also as in the proof of Proposition \ref{prop:tightemp} it can be argued that $B^1$ is a standard Brownian motion and thus we can conclude as in the proof of \eqref{eq:eq637}
\begin{equation}\label{eq:eq637b}
\tilde Z(t)\equiv (\tilde Q^1,\tilde Y^1,\tilde R^1)(t)=\Gamma\left(x+\int_0^\cdot
b_{\clr}(t,\tilde Q^1(t), \beta_{\clr}(t))dt +\sigma B^1(\cdot) \right)(t),\; t \in [0,T].
\end{equation}
Once again, by a standard martingale argument, one can argue that for $0 \le s \le t \le T$, $B^{1}(t) - B^1(s)$ is independent of
$$\tilde \clf_s \doteq \sigma \left\{ \tilde Q^1(s'), \tilde Y^1(s'), \tilde R^1(s'), \bar\beta_{\clr}(A \times [0,s']): 0\le s'\le s,\, A \in \clb(U)\right\}.$$
Thus denoting by $(\tilde \Omega, \tilde \clf, \tilde \PP)$ the probability space on which the limit processes are defined,
$\Xi = (\tilde \Omega, \tilde \clf, \{\tilde \clf_t\},\tilde \PP, B^1)$ is a system and
$(\bar\beta_{\clr}, \tilde Z) \in \cla_{\clr}(\Xi, 0, x, \bar \nu)$.
Exactly as in the proof of Proposition \ref{prop2} we see that the convergence in \eqref{eq:eq908} holds. Also using the weak convergence
$$(\tilde \nu^n, \tilde Q^{n,1}, \bar \beta^n_{\clr}) \Rightarrow (\bar \nu, \tilde Q^1, \bar \beta_{\clr})$$
and the Lipschitz property of $f$ in Assumption \ref{assumption1}, we have
\begin{align}
&\E\left[\int_{U\times [0,T]} f(t,\tnun(t), \tilde Q^{n,1}(t), u)\bar \beta^n_{\clr}(du\, dt) + g(\tnun(T),\tilde Q^{n,1}(T))\right]\nonumber\\
	&\quad \to \E\left[\int_0^T f(t,\bar \nu(t), \tilde Q^1(t), u)\bar \beta_{\clr}(du\, dt) + g(\bar \nu(T),\tilde Q^1(T))\right].\label{eq:eq504}
\end{align}
Combining the above convergence properties, we have as $n \to \infty$,
$$J^{n,1}(\tilde Q^n(0), \tilde \alpha^n_{-1}) \to J_{\bar \nu, \clr}(0, x, \bar \beta_{\clr}, \tilde Z) \ge V_{\bar \nu, \clr}(0,x) = V_{\bar \nu}(0,x).$$
Since the above holds for an arbitrary convergent subsequence of processes in \eqref{eq:eq426} and the sequence $\{\bar \beta^n_{\clr}\}$, the result follows. \hfill \qed

\skp
\noi
{\bf Acknowledgment.}
We are grateful to two anonymous referees for their suggestions that improved the presentation of the paper. We will also like to thank the associate editor  for the careful reading of the manuscript, valuable comments, and pointing out a key  error  in the convergence proof of the numerical scheme presented in a previous version of the paper. The corrected version of the statement and proof of the  convergence result is given in \cite{BBC2017}.

\small{
\bibliographystyle{abbrv} 
\bibliography{bib_Asaf_IMS} 
}\footnotesize{{\sc
\bigskip

\noindent
Erhan Bayraktar\\
Department of Mathematics\\
University of Michigan\\
Ann Arbor, MI 48109, USA\\
email: erhan@umich.edu\\
web: www-personal.umich.edu/$\sim$ erhan/
\skp

\noindent
Amarjit Budhiraja \\ 
Department of Statistics and Operations Research\\
University of North Carolina\\
Chapel Hill, NC 27599, USA\\
email: budhiraj@email.unc.edu
web: http://www.unc.edu/$\sim$ budhiraj/
\skp

\noindent
Asaf Cohen\\
Department of Statistics\\
University of Haifa\\
Haifa 31905, Israel\\
email: shloshim@gmail.com\\
web: https://sites.google.com/site/asafcohentau/

}}
\end{document}